\documentclass[journal, twoside]{IEEEtran}

\usepackage{setspace}
\usepackage{relsize}
\usepackage{url}

\usepackage{amsfonts}
\usepackage[cmex10]{amsmath} 
\usepackage{amssymb}

\usepackage{amsthm}	

\usepackage{mathrsfs}	
\usepackage[pdftex]{graphicx}	
\usepackage{epstopdf}	
\usepackage[breaklinks=true,linkcolor=blue,anchorcolor=black,citecolor=blue,urlcolor=blue]{hyperref}
\usepackage{color}
\usepackage{authblk}
\usepackage{cite}
\usepackage{multirow}
\usepackage[font=footnotesize]{caption}
\usepackage[font=footnotesize,subrefformat=parens]{subcaption}
\usepackage{hyperref}
\usepackage{array}									
\usepackage{wasysym}

\theoremstyle{plain} 
\newtheorem{theorem}{Theorem}
\newtheorem{proposition}{Proposition}
\newtheorem{corollary}{Corollary}

\theoremstyle{definition} 
\newtheorem{lemma}{Lemma}
\newtheorem{definition}{Definition}
\newtheorem{problem}{Problem}
\newtheorem{property}{Property}
\newtheorem{conjecture}{Conjecture}
\newtheorem{example}{Example}
\newtheorem{assumption}{Assumption}

\theoremstyle{remark} 
\newtheorem{remark}{Remark}
\newtheorem{fact}{Fact}
\newtheorem{claim}{Claim}

\newcommand{\theoremAlt}[2][]{		\begin{theorem}{{\textit{#1}}} 	{#2} \end{theorem}}
\newcommand{\corollaryAlt}[2][]{		\begin{corollary}{\textit{{#1}}}	{#2} \end{corollary}}

\newcommand{\assumptionAlt}[2][]{	\begin{assumption}{{#1}}		{#2} \end{assumption}}

\def\nn{\nonumber}
\def\beq{\begin{equation}}
\def\eeq{\end{equation}}
\def\beqs{\begin{equation*}}
\def\eeqs{\end{equation*}}
\def\balign{\begin{align}}
\def\ealign{\end{align}}
\def\bea{\begin{eqnarray}}
\def\eea{\end{eqnarray}}
\def\ba{\begin{array}}
\def\ea{\end{array}}
\def\bcenter{\begin{center}}
\def\ecenter{\end{center}}
\def\bitem{\begin{itemize}}
\def\eitem{\end{itemize}}
\def\benum{\begin{enumerate}}
\def\eenum{\end{enumerate}}
\def\bdoc{\begin{document}}
\def\edoc{\end{document}}
\def\defeq{{\stackrel{\Delta}{=}}}
\def\nn{\nonumber}

\newcommand{\eg}{\mbox{\textit{e.g.}, \/}}
\newcommand{\etal}{\mbox{\textit{et al.} \/}}
\newcommand{\viz}{\mbox{\textit{\iet viz.},\ \/}}
\newcommand{\ie}{\mbox{\textit{i.e.},\ \/}}
\newcommand{\etc}{\mbox{\textit{etc.},\ \/}}
\newcommand{\cf}{\mbox{\textit{cf.},\ \/}}
\newcommand{\vs}{\mbox{\textit{vs.}\ \/}}

\def\red{\textcolor{red}}
\def\blue{\textcolor{blue}}
\def\black{\textcolor{black}}
\def\yellow{\textcolor{yellow}}
\def\green{\textcolor{green}}
\def\tcbg{\textcolor{bgrd}}
\def\magenta{\textcolor{magenta}}
\def\white{\textcolor{white}}
\def\gray{\textcolor{gray}}
\newcommand\Quote[1]{\lq\textsl{#1}\rq}
\newcommand\fr[2]{{\textstyle\frac{#1}{#2}}}

\newcommand{\Ambb}{\mathbb{A}}
\newcommand{\Bmbb}{\mathbb{B}}
\newcommand{\Cmbb}{\mathbb{C}}
\newcommand{\Dmbb}{\mathbb{D}}
\newcommand{\Embb}{\mathbb{E}}
\newcommand{\Fmbb}{\mathbb{F}}
\newcommand{\Gmbb}{\mathbb{G}}
\newcommand{\Hmbb}{\mathbb{H}}
\newcommand{\Imbb}{\mathbb{I}}
\newcommand{\Jmbb}{\mathbb{J}}
\newcommand{\Kmbb}{\mathbb{K}}
\newcommand{\Lmbb}{\mathbb{L}}
\newcommand{\Mmbb}{\mathbb{M}}
\newcommand{\Nmbb}{\mathbb{N}}
\newcommand{\Ombb}{\mathbb{O}}
\newcommand{\Pmbb}{\mathbb{P}}
\newcommand{\Qmbb}{\mathbb{Q}}
\newcommand{\Rmbb}{\mathbb{R}}
\newcommand{\Smbb}{\mathbb{S}}
\newcommand{\Tmbb}{\mathbb{T}}
\newcommand{\Umbb}{\mathbb{U}}
\newcommand{\Vmbb}{\mathbb{V}}
\newcommand{\Wmbb}{\mathbb{W}}
\newcommand{\Xmbb}{\mathbb{X}}
\newcommand{\Ymbb}{\mathbb{Y}}
\newcommand{\Zmbb}{\mathbb{Z}}

\newcommand{\Amsc}{\mathcal{A}}
\newcommand{\Bmsc}{\mathcal{B}}
\newcommand{\Cmsc}{\mathcal{C}}
\newcommand{\Dmsc}{\mathcal{D}}
\newcommand{\Emsc}{\mathcal{E}}
\newcommand{\Fmsc}{\mathcal{F}}
\newcommand{\Gmsc}{\mathcal{G}}
\newcommand{\Hmsc}{\mathcal{H}}
\newcommand{\Imsc}{\mathcal{I}}
\newcommand{\Jmsc}{\mathcal{J}}
\newcommand{\Kmsc}{\mathcal{K}}
\newcommand{\Lmsc}{\mathcal{L}}
\newcommand{\Mmsc}{\mathcal{M}}
\newcommand{\Nmsc}{\mathcal{N}}
\newcommand{\Omsc}{\mathcal{O}}
\newcommand{\Pmsc}{\mathcal{P}}
\newcommand{\Qmsc}{\mathcal{Q}}
\newcommand{\Rmsc}{\mathcal{R}}
\newcommand{\Smsc}{\mathcal{S}}
\newcommand{\Tmsc}{\mathcal{T}}
\newcommand{\Umsc}{\mathcal{U}}
\newcommand{\Vmsc}{\mathcal{V}}
\newcommand{\Wmsc}{\mathcal{W}}
\newcommand{\Xmsc}{\mathcal{X}}
\newcommand{\Ymsc}{\mathcal{Y}}
\newcommand{\Zmsc}{\mathcal{Z}}

\newcommand{\Ascr}{\mathscr{A}}
\newcommand{\Bscr}{\mathscr{B}}
\newcommand{\Cscr}{\mathscr{C}}
\newcommand{\Dscr}{\mathscr{D}}
\newcommand{\Escr}{\mathscr{E}}
\newcommand{\Fscr}{\mathscr{F}}
\newcommand{\Gscr}{\mathscr{G}}
\newcommand{\Hscr}{\mathscr{H}}
\newcommand{\Iscr}{\mathscr{I}}
\newcommand{\Jscr}{\mathscr{J}}
\newcommand{\Kscr}{\mathscr{K}}
\newcommand{\Lscr}{\mathscr{L}}
\newcommand{\Mscr}{\mathscr{M}}
\newcommand{\Nscr}{\mathscr{N}}
\newcommand{\Oscr}{\mathscr{O}}
\newcommand{\Pscr}{\mathscr{P}}
\newcommand{\Qscr}{\mathscr{Q}}
\newcommand{\Rscr}{\mathscr{R}}
\newcommand{\Sscr}{\mathscr{S}}
\newcommand{\Tscr}{\mathscr{T}}
\newcommand{\Uscr}{\mathscr{U}}
\newcommand{\Vscr}{\mathscr{V}}
\newcommand{\Wscr}{\mathscr{W}}
\newcommand{\Xscr}{\mathscr{X}}
\newcommand{\Yscr}{\mathscr{Y}}
\newcommand{\Zscr}{\mathscr{Z}}

\def\Atiny{\mbox{\tiny{A}}}
\def\Btiny{\mbox{\tiny{B}}}
\def\Ctiny{\mbox{\tiny{C}}}
\def\Dtiny{\mbox{\tiny{D}}}
\def\Etiny{\mbox{\tiny{E}}}
\def\Ctiny{\mbox{\tiny{C}}}
\def\Ftiny{\mbox{\tiny{F}}}
\def\Gtiny{\mbox{\tiny{G}}}
\def\Htiny{\mbox{\tiny{H}}}
\def\Itiny{\mbox{\tiny{I}}}
\def\Jtiny{\mbox{\tiny{J}}}
\def\Ktiny{\mbox{\tiny{K}}}
\def\Ltiny{\mbox{\tiny{L}}}
\def\Mtiny{\mbox{\tiny{M}}}
\def\Ntiny{\mbox{\tiny{N}}}
\def\Otiny{\mbox{\tiny{O}}}
\def\Ptiny{\mbox{\tiny{P}}}
\def\Qtiny{\mbox{\tiny{Q}}}
\def\Rtiny{\mbox{\tiny{R}}}
\def\Stiny{\mbox{\tiny{S}}}
\def\Ttiny{\mbox{\tiny{T}}}
\def\Utiny{\mbox{\tiny{U}}}
\def\Vtiny{\mbox{\tiny{V}}}
\def\Wtiny{\mbox{\tiny{W}}}
\def\Xtiny{\mbox{\tiny{X}}}
\def\Ytiny{\mbox{\tiny{Y}}}
\def\Ztiny{\mbox{\tiny{Z}}}

\def\alphabf{\hbox{\boldmath$\alpha$\unboldmath}}
\def\betabf{\hbox{\boldmath$\beta$\unboldmath}}
\def\gammabf{\hbox{\boldmath$\gamma$\unboldmath}}
\def\deltabf{\hbox{\boldmath$\delta$\unboldmath}}
\def\epsilonbf{\hbox{\boldmath$\epsilon$\unboldmath}}
\def\zetabf{\hbox{\boldmath$\zeta$\unboldmath}}
\def\etabf{\hbox{\boldmath$\eta$\unboldmath}}
\def\iotabf{\hbox{\boldmath$\iota$\unboldmath}}
\def\kappabf{\hbox{\boldmath$\kappa$\unboldmath}}
\def\lambdabf{\hbox{\boldmath$\lambda$\unboldmath}}
\def\mubf{\hbox{\boldmath$\mu$\unboldmath}}
\def\nubf{\hbox{\boldmath$\nu$\unboldmath}}
\def\xibf{\hbox{\boldmath$\xi$\unboldmath}}
\def\pibf{\hbox{\boldmath$\pi$\unboldmath}}
\def\rhobf{\hbox{\boldmath$\rho$\unboldmath}}
\def\sigmabf{\hbox{\boldmath$\sigma$\unboldmath}}
\def\taubf{\hbox{\boldmath$\tau$\unboldmath}}
\def\upsilonbf{\hbox{\boldmath$\upsilon$\unboldmath}}
\def\phibf{\hbox{\boldmath$\phi$\unboldmath}}
\def\chibf{\hbox{\boldmath$\chi$\unboldmath}}
\def\psibf{\hbox{\boldmath$\psi$\unboldmath}}
\def\omegabf{\hbox{\boldmath$\omega$\unboldmath}}
\def\Sigmabf{\hbox{$\bf \Sigma$}}
\def\Upsilonbf{\hbox{$\bf \Upsilon$}}
\def\Omegabf{\hbox{$\bf \Omega$}}
\def\Deltabf{\hbox{$\bf \Delta$}}
\def\Gammabf{\hbox{$\bf \Gamma$}}
\def\Thetabf{\hbox{$\bf \Theta$}}
\def\Lambdabf{\mbox{$ \bf \Lambda $}}
\def\Xibf{\hbox{\bf$\Xi$}}
\def\Pibf{{\bf \Pi}}
\def\thetabf{{\mbox{\boldmath$\theta$\unboldmath}}}
\def\Upsilonbf{\hbox{\boldmath$\Upsilon$\unboldmath}}
\newcommand{\Phibf}{\mbox{${\bf \Phi}$}}
\newcommand{\Psibf}{\mbox{${\bf \Psi}$}}

\def\abf{{\bf a}}
\def\bbf{{\bf b}}
\def\cbf{{\bf c}}
\def\dbf{{\bf d}}
\def\ebf{{\bf e}}
\def\fbf{{\bf f}}
\def\gbf{{\bf g}}
\def\hbf{{\bf h}}
\def\ibf{{\bf i}}
\def\jbf{{\bf j}}
\def\kbf{{\bf k}}
\def\lbf{{\bf l}}
\def\mbf{{\bf m}}
\def\nbf{{\bf n}}
\def\obf{{\bf o}}
\def\pbf{{\bf p}}
\def\qbf{{\bf q}}
\def\rbf{{\bf r}}
\def\sbf{{\bf s}}
\def\tbf{{\bf t}}
\def\ubf{{\bf u}}
\def\vbf{{\bf v}}
\def\wbf{{\bf w}}
\def\xbf{{\bf x}}
\def\ybf{{\bf y}}
\def\zbf{{\bf z}}
\def\rbf{{\bf r}}
\def\xbf{{\bf x}}
\def\ybf{{\bf y}}
\def\Abf{{\bf A}}
\def\Bbf{{\bf B}}
\def\Cbf{{\bf C}}
\def\Dbf{{\bf D}}
\def\Ebf{{\bf E}}
\def\Fbf{{\bf F}}
\def\Gbf{{\bf G}}
\def\Hbf{{\bf H}}
\def\Ibf{{\bf I}}
\def\Jbf{{\bf J}}
\def\Kbf{{\bf K}}
\def\Lbf{{\bf L}}
\def\Mbf{{\bf M}}
\def\Nbf{{\bf N}}
\def\Obf{{\bf O}}
\def\Pbf{{\bf P}}
\def\Qbf{{\bf Q}}
\def\Rbf{{\bf R}}
\def\Sbf{{\bf S}}
\def\Tbf{{\bf T}}
\def\Ubf{{\bf U}}
\def\Vbf{{\bf V}}
\def\Wbf{{\bf W}}
\def\Xbf{{\bf X}}
\def\Ybf{{\bf Y}}
\def\Zbf{{\bf Z}}
\def\Ac{{\cal A}}
\def\Bc{{\cal B}}
\def\Cc{{\cal C}}
\def\Dc{{\cal D}}
\def\Ec{{\cal E}}
\def\Fc{{\cal F}}
\def\Gc{{\cal G}}
\def\Hc{{\cal H}}
\def\Ic{{\cal I}}
\def\Jc{{\cal J}}
\def\Kc{{\cal K}}
\def\Lc{{\cal L}}
\def\Mc{{\cal M}}
\def\Nc{{\cal N}}
\def\Oc{{\cal O}}
\def\Pc{{\cal P}}
\def\Qc{{\cal Q}}
\def\Rc{{\cal R}}
\def\Sc{{\cal S}}
\def\Tc{{\cal T}}
\def\Uc{{\cal U}}
\def\Vc{{\cal V}}
\def\Wc{{\cal W}}
\def\Xc{{\cal X}}
\def\Yc{{\cal Y}}
\def\Zc{{\cal Z}}

\def\0bf{\mathbf{0}}
\def\1bf{\mathbf{1}}

\def\lambdabf{\boldsymbol{\lambda}}
\def\mubf{\boldsymbol{\mu}}
\def\gammabf{\boldsymbol{\gamma}}
\def\xibf{\boldsymbol{\xi}}
\def\sigmabf{\boldsymbol{\sigma}}
\def\zetabf{\boldsymbol{\zeta}}

\def\Lambdabf{\boldsymbol{\Lambda}}
\def\Mubf{\boldsymbol{\Mu}}
\def\Gammabf{\boldsymbol{\Gamma}}
\def\Xibf{\boldsymbol{\Xi}}
\def\Sigmabf{\boldsymbol{\Sigma}}
\def\Zetabf{\boldsymbol{\Zeta}}

\def\1bf{\mathbf{1}}

\def\eqdef{\overset{\text{{\tiny def}}}{=}}

\DeclareMathOperator*{\cart}{\times}
\DeclareMathOperator*{\cov}{\Cmbb\textup{ov}}
\DeclareMathOperator*{\conv}{\text{conv}}
\DeclareMathOperator*{\diag}{\text{diag}}
\DeclareMathOperator*{\tr}{\textup{Tr}}
\DeclareMathOperator*{\minimize}{\text{minimize}}
\DeclareMathOperator*{\proj}{\text{proj}}
\DeclareMathOperator*{\toeplitz}{\text{toeplitz}}
\DeclareMathOperator*{\corr}{\text{corr}}
\DeclareMathOperator*{\rows}{\text{rows}}

\let\vec\relax
\DeclareMathOperator*{\vec}{\text{vec}}

\newcommand{\Top}{{\mbox{{\tiny $\top$}}}}

\renewcommand{\b}[1]{\ensuremath{\boldsymbol{\mathrm{#1}}}}

\def\sw{\textup{\textsf{sw}}}
\def\rp{\textup{\textsf{rs}}}
\def\cs{\textup{\textsf{cs}}}
\def\so{\textup{\textsf{so}}}
\def\vp{\textup{\textsf{vp}}}

\def\figwid{4cm} 
\newcommand*{\MyLocalPath}{}%

\newcommand{\ifThesis}[2]{\ifdefined\IEEEPARstart%
#2%
\else%
#1%
\fi}

\renewcommand{\qedsymbol}{\blacksquare}
\newcommand{\qedadhoc}{\tag*{$\blacksquare$}}

\providecommand{\boldsymbol}[1]{\mbox{\boldmath $#1$}}

\newenvironment{proofof}[1]{\bigskip \noindent {\bf Proof of #1:}\ }{\qed}

\IEEEoverridecommandlockouts

\begin{document}

\title{On the Efficiency of Connection Charges---\\Part II: Integration of Distributed Energy Resources}

\author{Daniel~Munoz-Alvarez,~Juan~F.~Garcia-Franco
		and~Lang~Tong
\thanks{This work is supported in part by the National Science Foundation under Grants CNS-1135844 and 15499.}
\thanks{D. Munoz-Alvarez and L. Tong are with the School of Electrical and Computer Engineering, Cornell University, Ithaca, NY, 14853, USA.  Emails:  {\tt\small dm634@cornell.edu}, {\tt\small lt35@cornell.edu}}
}

\markboth{}%
{}
%

\maketitle 

\begin{abstract}
This two-part paper addresses the design of retail electricity tariffs for distribution systems with distributed energy resources (DERs).
Part I presents a framework to optimize an ex-ante two-part tariff for a regulated monopolistic retailer who faces stochastic wholesale prices on the one hand and stochastic demand on the other.
In Part II, the integration of DERs is addressed by analyzing their endogenous effect on the optimal two-part tariff and the induced welfare gains.
Two DER integration models are considered: $(i)$ a decentralized model involving behind-the-meter DERs in a net metering setting, and
$(ii)$ a centralized model involving DERs integrated by the retailer.

It is shown that DERs integrated under either model can achieve the same social welfare and the net-metering tariff structure is optimal.
The retail prices under both integration models are equal and reflect the expected wholesale prices.
The connection charges differ and are affected by the retailer's fixed costs as well as the statistical dependencies between wholesale prices and behind-the-meter DERs.
In particular, the connection charge of the decentralized model is generally higher than that of the centralized model.

An empirical analysis is presented to estimate the impact of DER on welfare distribution and inter-class cross-subsidies using real price and demand data and simulations.
The analysis shows that, with the prevailing retail pricing and net-metering, consumer welfare decreases with the level of DER integration.
Issues of cross-subsidy and practical drawbacks of decentralized integration are also discussed.

\end{abstract}

\begin{IEEEkeywords}
Retail tariff design, dynamic pricing, connection charges, distributed energy resources, renewables, storage.
\end{IEEEkeywords}

\section{Introduction} \label{sec:intro2}

\ifdefined\IEEEPARstart
	\IEEEPARstart{T}{his}
\else
	This
\fi
two-part paper studies the design of dynamic retail electricity tariffs for distribution systems with distributed renewable and storage resources.
We consider a regulated monopolistic retailer who, on the one hand, serves residential customers with stochastic demands, and on the other hand, interfaces with an exogenous wholesale market with stochastic prices.
In this framework, we analyze both customer-integrated and retailer-integrated distributed energy resources (DERs).
Our goal is to shed lights on the widely adopted net metering compensation mechanism and the efficiency loss implied by some of the prevailing retail tariffs when an increasing amount of DERs are integrated into the distribution system.

While Part I \cite{MunozTong16partIarxiv} establishes a framework to analyze the efficiency of revenue adequate tariffs with connection charges, Part II 
extends it to address the integration of DERs.

The main contribution of Part II is twofold.
First, we characterize analytically the optimal revenue adequate ex-ante two-part tariff for a distribution system with renewables and storage integrated by customers or the retailer.
We characterize the consumer (and social) welfare achieved by the optimal two-part tariff under both integration models.
This analysis is an application of the classical Ramsey pricing theory \cite{BrownSibley86} with extensions to accommodate the multi-period integration of stochastic DERs.
Second, we analyze a numerical case study based on empirical data that estimates the increasingly larger inefficiencies and interclass cross-subsidies caused by DERs when net metering tariffs with price markups are used to maintain revenue adequacy.
In this context, the derived optimal two-part tariffs and a centralized DER integration model offer two alternatives to mitigate these undesirable effects.

The main results of Part II are as follows.
We leverage the retail tariff design framework established in \cite{MunozTong16partIarxiv} to accommodate the integration of DERs by customers (in a net-metering setting) and by the retailer in Section \ref{sec:retailTariffDesignDER}.
The extended framework considers heterogeneous customers with arbitrary behind-the-meter renewables and storage.
Therein, we derive the optimal ex-ante two-part tariff under both DER integration models and the combined effect of this tariff and DERs on consumer and social welfare.

We find that under the optimal two-part tariff, DERs integrated under either model bring the same gains in social and consumer welfare.
This is in contrast to prevailing volumetric tariffs under which the integration of DERs can increase or decrease social and consumer welfare depending critically on the integration model and the retailer's fixed costs.
Indeed, we demonstrate that the two-part tariff structure is optimal in the sense that no other tariff structure ---however complex--- can achieve a strictly higher social welfare.
This means that the two-part \emph{net metering} tariff of the decentralized model is optimal as a DER compensation mechanism.

These welfare effects are explained by the structure of the optimal ex-ante two-part tariff.
We show that under both integration models the derived tariff consists of an identical time-varying price and a distinct connection charge.
In particular, the time-varying price reflects the wholesale prices and their statistical correlation with the elasticity of the random demand.
The optimal connection charge allocates uniformly among customers the retailer's fixed costs and additional costs and savings caused by risks and the integrated DERs.
Indeed, while savings from retailer-integrated DERs reduce the connection charge, customer-integrated DERs induce slight increments or reductions caused by risks introduced by renewables.

The theoretical analysis of DER integration is complemented in Section \ref{sec:caseStudy2} with an empirical study based on publicly available data from NYISO and the largest utility company in New York City.
The performance of the optimal ex-ante two-part tariffs is compared with several other ex-ante tariffs for different levels of DER penetration, under both integration models.
Tariffs used as benchmarks include the optimal linear tariff and two-part flat tariffs used extensively in practice by utilities.
In particular, relative to a base case with a nominal two-part flat tariff and no DERs, we estimate the efficiency gains or losses brought by tariff changes in Section \ref{sec:caseStudy2:baseCase}.
Subsequently, in Section \ref{sec:caseStudy2:efficiency}, we estimate the efficiency gains or loses brought by the integration of DERs under both integration models and the various ex-ante considered tariffs.
Most notably, our results estimate that the efficiency gains brought by switching from flat to hourly pricing, which are below $1\%$ (of the utility's gross revenue) for most relevant cases, can be more than tripled by a $\$10$ increase in the monthly connection charge.
Moreover, for the case with customer-integrated DERs, we estimate in Section \ref{sec:caseStudy2:xsubsidies} the indirect cross-subsidies that customers without DERs give to DER-owning customers due to net metering tariffs with marked-up retail prices.
All our estimations in this case study assume a stylized model for thermostatically controlled loads.

Concluding remarks and proof sketches of the main results are included in Section \ref{sec:conclusions2} and the Appendix, respectively.
Detailed proofs of all results can be found in \cite{MunozTong16partIIarxiv}.

\

\vspace{-10pt}

\subsection{Related Work}

The literature on retail electricity tariff design is extensive \cite{Munoz16}, and there is an increasing interest in addressing the integration of DERs.
We briefly discuss works that are relevant to our paper.
Based on their main focus, we group these works into two categories: $(i)$ tariff design for fixed cost recovery with DERs, and $(ii)$ optimal demand response with DERs.

\subsubsection{Tariff design for fixed cost recovery with DERs}
The general principles used in retail tariff design are briefly reviewed in \ifThesis{\cite{RodriguezEtal08, RenesesRodriguez14}}{\cite{RenesesRodriguez14}} and more extensively in \ifThesis{\cite{BraithwaitEtal07, LBNL16b}}{\cite{LBNL16b}}, and the additional challenges brought by DERs are discussed in \cite{Costello15}.
In the light of such challenges, current tariff design practices and broader regulatory issues are being revised in comprehensive studies to address the adoption of DERs%
\ifThesis{
\cite{NREL13,MIT15,DPS15_Track2WhitePaper,NREL15,LBNL16a,LBNL16b,NARUC16}, to estimate the impact of different tariff structures on the bills of residential customers with solar PV \cite{NREL15}, and to investigate pricing issues related to the interaction between distribution utilities and the owners of DERs \cite{LBNL16a}.
}{\cite{NREL13,MIT15,NARUC16}}. 

Research efforts to study more specific issues of DER integration such as \cite{EidEtal14,JargstorfBelmans15,DarghouthEtal16,Sioshansi16} have also emerged. 
For instance, in \cite{JargstorfBelmans15}, the trade-off between multiple tariff design criteria is studied in a multi-objective optimization framework.
An analytical approach leverages a generation capacity investment model in \cite{Sioshansi16} to characterize sufficient conditions for RTP and flat tariffs to be revenue adequate.
More empirical approaches are conducted in \cite{EidEtal14}, where interclass cross-subsidies and revenue shortfalls caused by net metering tariffs are estimated, and in \cite{DarghouthEtal16}, which estimates the impact of tariff structure and net metering on the deployment of distributed solar PV.

Finally, there is an increasing volume of literature studying the ``death spiral'' of DER adoption \cite{ChewEtal12, CaiEtal13, Kind13, MIT13, CostelloHemphill14, RMI14, RMI15}, which is presented as a threat on the financial viability of utilities.
This threat refers to a self-reinforcing feedback loop of DER adoption involving a decline in energy sales
and the persistent attempt to recover utilities' fixed costs by increasing volumetric charges.
The empirical analysis in \cite{CaiEtal13}, for example, models the effect that price feedback loops may have on the adoption of solar PV and concludes that it may not be significant within the next decade.
In \cite{MIT13}, an extensive list of factors that affect the system dynamics of DER adoption is presented.
It concludes that while the feedback loop is possible, it is not predetermined and can be avoided.
A stylized demand model is used in \cite{CostelloHemphill14} to argue that a minimum of price elasticity is required for the threat to be an actual problem.
The work in \cite{RMI15} provides an estimate of the evolution of the lowest-cost configuration (namely grid only, grid+solar, or grid+solar+battery) for residential and commercial customers to satisfy their load in the long-term for a few U.S. cities.

There are still important gaps in this subject.
For example, none of the works above studies the efficiency loss and, with the exception of \cite{EidEtal14}, the interclass cross-subsidies entailed by the adoption of DERs under net metering tariffs.
This is precisely a focus of our work.  

\subsubsection{Optimal demand response with DERs}
Many works focus on deriving optimal retail pricing schemes to induce desired electricity consumption behavior on customers with DERs such as \cite{ChenEtal12, TangEtal14, JiaTong16b, HanifEtal16}. 
For example, in \cite{JiaTong16b}, the authors consider customer and retailer integrated renewables and storage separately in a setting similar to ours.
They derive dynamic linear tariffs that maximize an objective that balances the retailer profit and customers' welfare.
Unlike our work, however, none of these works consider explicitly a revenue adequacy constraint nor the use of connection charges.

\section{Retail Tariff Design with DERs} \label{sec:retailTariffDesignDER}


%
%
%
%


\subsection{Multi-period Ramsey Pricing under Uncertainty}

Consider a regulator who sets a retail electricity tariff $T$ in advance (ex-ante) to maximize the welfare of $M$ customers over a billing cycle of $N$ time periods, subject to a net revenue sufficiency constraint for the monopolistic retailer serving the load.
Expectations are used to deal with the uncertainties that naturally arise when fixing a tariff in advance of actual usage.

To quantify the customers' welfare we use the notion of consumers' surplus, which measures the difference between the gross benefit derived from consumption and what the customer pays for it.
Formally, we assume that given a tariff $T$, customer $i$ consumes a profile $q^i(T,\omega^i) \in \Rmbb^N$ within the $N$-period billing cycle contingent on the random evolution of the local state $\omega^i = (\omega^i_1,\ldots,\omega^i_N) \in \Rmbb^N$, provided that $q^i$ is purchased from the retailer.
Accordingly, customer $i$ derives an expected surplus
\begin{align} \label{eq2:cs}
\overline{\cs}^i(T) = \Embb \big[ S^i(q^i(T,\omega^i),\omega^i) - T(q^i(T,\omega^i)) \big], 
\end{align}
where $T:\Rmbb^N \rightarrow \Rmbb$ and $S^i(q^i(T,\omega^i),\omega^i)$ is the derived gross benefit.
Collectively, customers derive an expected consumer surplus $\overline{\cs}(T) = \Embb[ \sum_{i=1}^M \cs^i(T) ]$, where the expectation is taken with respect to the $M$-tuple $\omega = (\omega^1,\ldots,\omega^M)$.

Similarly, the expected retailer surplus or net revenue is
\begin{align} \label{eq2:rs}
\overline{\rp}(T) = \Embb \big[ \mbox{$\sum_{i=1}^M$} T(q^i(T,\omega^i)) - \lambda^{\Top} q(T,\omega) \big], 
\end{align}
where $\lambda \in \Rmbb^N$ is the profile of random real-time wholesale prices, $q(T,\omega)$ is the aggregated demand profile, $\lambda^{\Top} q(T,\omega)$ is the energy cost faced by the retailer, and the expectation is over the uncertain evolution of the global state $\xi=(\lambda,\omega)$.

Adding the consumer and retailer surplus together yields the (expected) social surplus
$
\overline{\sw}(T) = \overline{\cs}(T) + \overline{\rp}(T)
$
which quantifies the social welfare induced by a tariff $T$.

We can now formulate the regulator's tariff design problem as the optimization problem
\begin{align} \label{eq2:reg problem}
\max_{T(\cdot)}	&	\	\overline{\cs}(T) \quad \text{s.t.} \quad \overline{\rp}(T) = F,
\end{align}
where $F$ is a constant representing the non-energy costs faced by the retailer that need to be passed on to its customers\footnote{$F$ may include delivery, metering, and customer service costs, and it may also recognize that a regulated firm should be allowed to earn some profit.}.
As such, \eqref{eq2:reg problem} is a version of the Ramsey pricing problem\footnote{Ramsey pricing is pricing efficiently subject to a breakeven constraint \cite{BrownSibley86}.
With \eqref{eq2:reg problem}, we seek to apply Ramsey pricing to a single service with time-varying, random marginal costs and temporally dependent stochastic demands.
}.
In particular, we consider ex-ante two-part tariffs\footnote{This restriction may involve no loss of generality (see Thm. \ref{thm:optimalityDERs} below).
} $T(q)=A+\pi^{\Top}q$ with connection charge $A \in \Rmbb$ and time-varying price $\pi \in\Rmbb^N$.
These tariffs induce an individual consumption profile $q^i(T,\omega^i) = D^i(\pi,\omega^i)$, where $D^i(\cdot,\omega^i)$ is a demand function assumed to be nonnegative, continuously differentiable in $\pi$, and with a negative definite Jacobian $\nabla_{\pi} D^i(\pi,\omega^i) \in \Rmbb^{N \times N}$ that satisfies the following assumption\footnote{A detailed discussion on the implications of this assumption and special cases when it is satisfied can be found in Part I \cite{MunozTong16partIarxiv}.}.
\assumptionAlt{ \label{Assumption 1:part2}
$g(\pi)=\Embb[\nabla_{\pi} D(\pi,\omega) (\pi-\lambda)]$ is such that the Jacobian matrix $\nabla g(\pi)$ is negative definite (nd).
}

In the following sections we accommodate the integration of DERs into the tariff design framework above.
To that end, we assume that either customers or the retailer have access to distributed renewable and storage resources.
We model an agent's access to renewables as the ability to use a state-contingent energy profile $r\in\Rmbb_+^N$ at no cost.
Similarly, we model access to a storage with capacity $\theta \in \Rmbb_+$ as the ability to offset energy \emph{needs} with any vector of storage discharges $s \in \Rmbb^N$ in the operation constraint set\footnote{The lossless storage model defined by $\Umsc(\theta)$, which assumes no initial charge nor charging/discharging rate limits, involves no loss of generality since more complex storage models can be accommodated redefining $\Umsc(\theta)$.}
$$
\Umsc(\theta) = \left\{ s \in \Rmbb^N \ \Big\vert \ 0 \leq - \mbox{$\sum_{t=1}^k$} s_t \leq \theta, \ k=1,\ldots,N \right\}.
$$
We define the (arbitrage) \emph{value} of the storage given a deterministic price vector $\pi \in\Rmbb^N$ as
\begin{align} \label{eq2:storage operation}
V^{\textsc{s}}(\pi, \theta) = \max_{s\in\Rmbb^N} \left\{ \pi^{\Top}s \ \Big\vert \ s \in \Umsc(\theta) \right\},
\end{align}
and let $s^{*}(\pi,\theta)$ denote an optimal solution of \eqref{eq2:storage operation}.

In what follows, we focus on characterizing solutions to problem \eqref{eq2:reg problem} considering DERs integrated either behind the meter by customers in a net-metering setting or by the retailer.

\subsection{Decentralized (behind-the-meter) DER Integration} \label{sec:decentralizedIntegration}

Suppose that customers install renewables and a battery behind the meter.
Let $r^i(\omega^i) \in \Rmbb^N$ and $s^i \in \Rmbb^N$ denote the energy customer $i$ obtains from renewable resources in state $\omega^i$ and from the battery, respectively, and let $\theta^i \in \Rmbb_+$ represent its storage capacity.
We operate in a net-metering setting where tariffs depend only on
$
d^i = q^i - r^i - s^i,
$
which we use to represent customer $i$'s net-metered demand.
Hence, given a tariff $T(d)=A+\pi^{\Top}d$, customer $i$ chooses consumption $q^i_k$ and storage operation $s^i_k$ at each time $k$ contingent on $\omega^i_1,\ldots,\omega^i_k$
to solve the multistage stochastic program
\begin{subequations} \label{eq2:cs with DERs}
\begin{align}
\overline{\cs}^i(T) = \max_{q^i(\cdot),s^i(\cdot)} &\ \Embb[S^i(q^i(\omega^i),\omega^i)-T(d^i(\omega^i))],  \label{eq2:cs with DERs:objective} \\
							\text{s.t} \ \ &\ s^i(\omega^i) \in \Umsc(\theta^i). \label{eq2:cs with DERs:cons}
\end{align}
\end{subequations}
A key observation is that the linearity of two-part tariffs implies that customer $i$'s problem \eqref{eq2:cs with DERs} can be separated into two sub-problems: choosing $q^i(\cdot)$ to maximize $\Embb[S^i(q^i(\omega^i),\omega^i)-\pi^{\Top} q^i(\omega^i)]-A$ and choosing $s^i(\cdot)$ to maximize $\Embb[\pi^{\Top} s^i(\omega^i)]$ subject to \eqref{eq2:cs with DERs:cons}.
The former problem is equivalent to that of customers without DERs analyzed in \cite{MunozTong16partIarxiv}, whose solution characterizes the demand function $D^i(\pi,\omega^i)$.
As for the second sub-problem, it is clear from \eqref{eq2:storage operation} that $s^{*}(\pi,\theta^i)$ is an optimal solution.
These solutions constitute an optimal solution to \eqref{eq2:cs with DERs} and thus a net demand function
\begin{align} \label{eq2:demand with DERs}
d^i(\pi,\omega^i) = D^i(\pi,\omega^i)-r^i(\omega^i)-s^{*}(\pi,\theta^i).
\end{align}
This fundamental separation of the customer's problem yields the following result, where we use $r(\omega) = \sum_{i=1}^M r^i(\omega^i)$ and $s = \sum_{i=1}^M s^i$ for notational convenience.

\theoremAlt{\label{thm:decentralized}
Suppose that customers have access to renewables and storage as characterized in \eqref{eq2:cs with DERs} and \eqref{eq2:demand with DERs}.
If $\nabla_{\pi} D(\pi,\omega)$ and $\lambda$ are uncorrelated\footnote{The absence of decentralized storage makes this condition unnecessary.}, then the two-part tariff $T^*_{\textsc{dec}}$ that solves problem \eqref{eq2:reg problem} is given by $\pi^*_{\textsc{dec}} = \overline{\lambda}$ and
\begin{align} \label{eq2:A decentralized}
A^*_{\textsc{dec}} 	&= A^* - \mbox{$\frac{1}{M}$} \tr\left( \cov \left( \lambda, r(\omega)  \right) \right),
\end{align}
where $A^*$, the connection charge in the absence of DER, would be given by
$
A^* = \mbox{$\frac{1}{M}$} \left( F + \tr(\cov(\lambda,D(\overline{\lambda},\omega))) \right).
$
}

%
%
%
%
%
%

Before discussing some implications of Theorem \ref{thm:decentralized}, we examine the condition that $\nabla_{\pi} D(\pi,\omega)$ and $\lambda$ are uncorrelated.
This condition holds in many situations.
In particular, it holds for demands that are not much affected by consumers' local randomness, such as the charging of electric vehicles and typical household appliances.
It even holds for smart HVAC loads that are affected by random temperature fluctuations since their demand takes the form $D(\pi,\omega)=D(\pi) + b(\omega)$, \ie a demand with additive disturbances \cite{JiaTong16a}.



The tariff $T^*_{\textsc{dec}}$ in Thm. \ref{thm:decentralized} reveals the following.
Letting retail prices reflect an unbiased estimate of the marginal costs of electricity ($\lambda$) maximizes social and consumer welfare.
Under net metering, this implies that the retailer should buy customers' energy surplus (from DERs) at the same price that he buys energy at the wholesale market (in expectation).

The expression for $A^*_{\textsc{dec}}$ in \eqref{eq2:A decentralized} has an intuitive interpretation.
It indicates that the integration of behind-the-meter DERs would require adjustments to the connection charge.
These adjustments could be positive if the integrated renewables tend to cause wholesale prices to drop (\ie negative correlation), but they could be negative otherwise.
Consequently, these adjustments can increase or decrease the consumer surplus of customers without DERs because the former are perceived by \emph{all} customers as changes in their electricity bills.

The welfare gains brought by decentralized DERs depend critically on retail tariffs.
To assess the performance of two-part tariffs in this regard we first need a point of comparison.
In the absence of DERs, $T^*_{\textsc{dec}}$ reduces to the optimal ex-ante two-part tariff $T^*(q)=A^*+\pi^{*\Top}q$ derived in \cite{MunozTong16partIarxiv}, where $\pi^*=\overline{\lambda}$ under the assumption in Theorem \ref{thm:decentralized}.
As a point of comparison, consider that in the absence of DERs and under tariff $T^*$, customers derive an expected surplus $\overline{\cs}_0(T^*)=\overline{\sw}^{*}_0 - F$, the retailer derives $\overline{\rp}_0(T^*)=F$, and social welfare is
\begin{align} \label{eq2:sw0 star}
\overline{\sw}^{*}_0 &= \sum_{i=1}^M \Embb \big[ S^i(D^i(\pi^*,\omega^i),\omega^i) - \lambda^{\Top}D^i(\pi^*,\omega^i) \big].
\end{align}

\corollaryAlt{\label{cor:decentralized}
Under the tariff $T^*_{\textsc{dec}}$, customer-integrated DERs induce an expected total surplus $\overline{\sw}(T^*_{\textsc{dec}}) = \overline{\sw}^*_0 + \sum_{i=1}^M V^{\textsc{s}}(\overline{\lambda},\theta^i) + \Embb[\lambda^{\Top} r^i(\omega^i)]$ that is independent of $F$ and 
$\overline{\cs}(T^*_{\textsc{dec}}) = \overline{\cs}_0(T^*) + \sum_{i=1}^M V^{\textsc{s}}(\overline{\lambda},\theta^i) + \Embb[\lambda^{\Top} r^i(\omega^i)]$.
}

The expressions $\overline{\cs}_0(T^*)$ and $\overline{\cs}(T^*_{\textsc{dec}})$ above characterize the tradeoff between the retailer's surplus target $F$ and consumers' surplus $\overline{\cs}$ induced by the tariffs $T^*$ and $T^*_{\textsc{dec}}$, respectively.
Indeed, noting the linear dependence of $A^*$ in $F$, it becomes clear that in both cases the $\overline{\rp}$-$\overline{\cs}$ tradeoff is linear, as illustrated in Fig. \ref{fig:PFdecentralized}.
Moreover, the fact that the social welfare achieved in both cases ($\overline{\sw}^{*}_0$ and $\overline{\sw}(T^*_{\textsc{dec}})$) does not depend on $F$ implies that said tradeoff is not only linear but one-to-one (\ie the Pareto fronts in Fig. \ref{fig:PFdecentralized} have slope $-1$).
This means that while an increased net revenue target $F+\Delta F$ decreases consumer surplus in expectation, it does not decrease social surplus.
Conversely, the integration of DERs behind-the-meter increases both social and consumer surplus by $\sum_{i=1}^M V^{\textsc{s}}(\overline{\lambda},\theta^i) + \Embb[\lambda^{\Top} r^i(\omega^i)]$ in expectation regardless of the retailer's net revenue target $F$.

\begin{figure}[t]
\centering
\subcaptionbox{Optimal two-part tariffs.\label{fig:example:CSgains:a}}
{\includegraphics[width=0.5\linewidth]{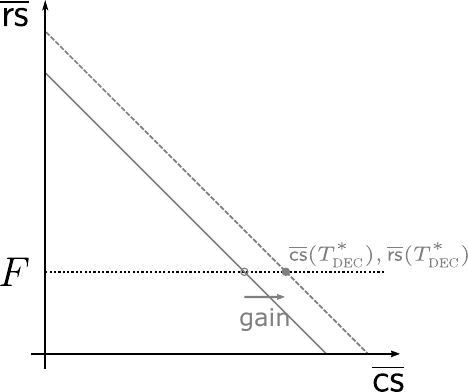}}
\caption{Efficient Pareto fronts or $\overline{\rp}$-$\overline{\cs}$ tradeoff induced by optimal ex-ante two-part tariffs $T^*_{\textsc{dec}}$ (solid line) and $T^*$ (dashed line) with and without DERs.}
\label{fig:PFdecentralized}
\end{figure}


Another implication of the optimal two-part tariff $T^*_{\textsc{dec}}$ is the likely impact it would have on the rapid adoption of behind-the-meter DERs.
Prevailing tariffs that rely on retail markups to achieve revenue adequacy provide an strong incentive for customers to integrate Distributed Generation (DG).
This is because, under net-metering, the higher the retail prices, the more savings DG represents.
By eliminating retail markups and imposing virtually unavoidable connection charges, $T^*_{\textsc{dec}}$ generally reduces such savings.
Hence, $T^*_{\textsc{dec}}$ is likely to decelerate the adoption of decentralized DERs compared to the prevailing less efficient retail tariffs.
This suggests that there is a tradeoff between efficiency and the rapid adoption of behind-the-meter DERs.

\subsubsection{Optimality of Net Metering} \label{sec:performance with DER}

We have restricted the regulator to offer net-metering two-part tariffs.
There are, however, alternative mechanisms to compensate DERs that do not rely on net metering (\eg feed-in tariffs).
We argue that the regulator cannot improve upon the efficiency attained by $T^*_{\textsc{dec}}$ with more complex ex-ante tariffs under certain condition\ifThesis{\footnote{The optimality argument in \cite[Sec. III.C]{MunozTong16partIarxiv} without DERs applies to the case with retailer-integrated DERs presented in the following section since both problems are equivalent except for a difference in the parameter $F$.}}{}.
This holds true because $T^*_{\textsc{dec}}$ induces the same efficiency attained by the social planner, which provides an upper bound to the regulator's problem \eqref{eq2:reg problem}\ifThesis{\footnote{Due to the restriction to ex-ante tariffs, we restrict the social planner's decisions to be contingent on each customer's local state $\omega^i$.
This is because ex-ante tariffs cannot carry updated information of the global state $\xi=(\lambda,\omega)$, unlike real-time or ex-post tariffs.}}{}.


\theoremAlt{ \label{thm:optimalityDERs}	
Suppose that customers have access to renewables and storage.
If wholesale prices $\lambda$ and customers' states $\omega$ are statistically independent (\ie $\lambda \perp \omega$), then $T^*_{\textsc{dec}}$ is an optimal solution of \eqref{eq2:reg problem} among the class of ex-ante tariffs.}

Namely, the restriction to two-part net metering tariffs, which are simple and thus practical tariffs, imply no loss of efficiency if $\lambda \perp \omega$.
The latter condition, however, makes the result somewhat restrictive as it applies to loads not affected by customers' local randomness such as washers and dryers, computers, batteries and EV charging \emph{but} not to HVAC loads or behind-the-meter solar and wind DG. 
Nonetheless, said condition suggests that if the net load and $\lambda$ are poorly correlated (or either exhibits little uncertainty at the time the tariff is fixed) then $T^*_{\textsc{dec}}$ may have a good performance.

\subsection{Centralized (retailer-based) DER Integration} \label{sec:centralizedIntegration}

As an alternative to behind-the-meter DERs, we now consider the case where the retailer installs DERs within the distribution network.
To that end, suppose that the retailer has access to a renewable supply $r^o(\xi)\in\Rmbb^N_+$ and a storage capacity $\theta^o \in \Rmbb_+$.
Without loss of generality, we assume that the retailer determines the operation of storage before the billing cycle starts (\ie ex-ante)\footnote{Allowing storage operation to be contingent on partial observations of $\lambda\in\Rmbb^N$ (say $s^o(\lambda)\in\Rmbb^N$) only makes the maximum value of $\Embb[\lambda^{\Top}s^o(\lambda)]$ over $s^o(\lambda)\in \Umsc(\theta)$ in \eqref{eq2:retailer separation} hard to compute under general assumptions.}.
Assuming that the retailer operates storage to maximize his net revenue, the resulting surplus induced by a tariff $T$ can be written as
\begin{small}
\begin{align}
\overline{\rp}(T)	&= \max_{s \in \Umsc(\theta^o)} \Embb\left[ \sum_{i=1}^{M} T(q^{i}(T,\omega^i)) \ - \lambda^{\Top} (q(T,\omega)-r^o(\omega)-s) \right]  \nn \\
						&= \overline{\rp}_0(T) - \Embb[\lambda^{\Top} r^o(\xi)] - V^{\textsc{s}}(\overline{\lambda},\theta^o),  \label{eq2:retailer separation}
\end{align}
\end{small}%
The fact that the two last terms in \eqref{eq2:retailer separation} do not depend on $T$ facilitates obtaining the following result since both terms simply offset the surplus target $F$ when imposing $\overline{\rp}(T) = F$.

\theoremAlt{\label{thm:centralized}
Suppose that the retailer has access to renewables and storage as characterized in \eqref{eq2:retailer separation}.
Then the two-part tariff $T^*_{\textsc{cen}}$ that solves problem \eqref{eq2:reg problem} is given by
\begin{align} 
\pi^*_{\textsc{cen}} = \overline{\lambda} + \Embb[\nabla_{\pi}D(\pi^*_{\textsc{cen}},\omega)]^{-1} \Embb[\nabla_{\pi}D(\pi^*_{\textsc{cen}},\omega) (\lambda-\overline{\lambda})], \nn
\end{align}
\begin{align} \label{eq2:A centralized}
A^*_{\textsc{cen}} = A^* - \mbox{$\frac{1}{M}$} \big( V^{\textsc{s}}(\overline{\lambda},\theta^o) + \Embb[\lambda^{\Top} r^o(\xi)] \big)
\end{align}
where $A^*$, the connection charge in the absence of DER, would be given by
$
A^* = \mbox{$\frac{1}{M}$} \left( F - \Embb[(\pi^*_{\textsc{cen}} - \lambda)^{\Top}D(\pi^*_{\textsc{cen}},\omega)] \right).
$
}

%
%
%
%
%
%
%

We first note that, unlike Thm. \ref{thm:decentralized}, Thm. \ref{thm:centralized} does not require $\nabla_{\pi} D(\pi,\omega)$ and $\lambda$ to be uncorrelated.
However, if this condition is satisfied it holds that $\pi^*_{\textsc{cen}}=\pi^*_{\textsc{dec}}=\overline{\lambda}$.
In other words, under optimally set ex-ante two-part tariffs, the integration of DERs by either customers or their retailer do not require updating prices to maintain revenue adequacy.
Hence, in both cases, any potential feedback loop of DER integration on retail prices (and thus on consumption) is undermined.


In terms of the connection charge in \eqref{eq2:A centralized}, the integration of DERs by the retailer results in reductions relative to $A^*$.
These reductions contrast with the potential increments required by customer-integrated DERs (\cf $A^*_{\textsc{dec}}$ in \eqref{eq2:A decentralized}).
The underlying reason for such difference is intuitive, specially considering the identical retail prices $\pi^*_{\textsc{cen}}=\pi^*_{\textsc{dec}}$.
Decentralized DERs represent savings in volumetric charges for customers (with reduced net loads) whereas centralized DERs represent savings in electricity purchases for the retailer.
Because the latter savings cannot increase the retailer surplus beyond $F$, they are allocated uniformly between customers through reductions in the connection charge.


Unlike with decentralized DERs in general, the welfare gains brought by DERs integrated (and operated) by the retailer do not depend on retail tariffs.
This is formalized by the following result.
\corollaryAlt{\label{cor:centralized}
Under the tariff $T^*_{\textsc{cen}}$, retailer-integrated DERs induce an expected total surplus $\overline{\sw}(T^*_{\textsc{cen}})= \overline{\sw}^*_0 + V^{\textsc{s}}(\overline{\lambda},\theta^o) +\Embb[\lambda^{\Top} r^o(\xi)]$ that is independent of $F$ and 
$\overline{\cs}(T^*_{\textsc{cen}}) = \overline{\cs}_0(T^*) + V^{\textsc{s}}(\overline{\lambda},\theta^o) +\Embb[\lambda^{\Top} r^o(\xi)]$.
}

In Cor. \ref{cor:centralized}, $\overline{\cs}(T^*_{\textsc{cen}})$ characterizes a linear one-to-one tradeoff between $F$ and $\overline{\cs}$ induced by the $T^*_{\textsc{cen}}$.
This tradeoff---equivalent to the Pareto front induced by $T^*_{\textsc{dec}}$ in Fig. \ref{fig:PFdecentralized}---is characterized by the social welfare achieved by $T^*_{\textsc{cen}}$, $\overline{\sw}(T^*_{\textsc{cen}})$.
Consequently, similar to behind-the-meter DERs, the integration of centralized DERs increases both social and consumer surplus by $V^{\textsc{s}}(\overline{\lambda},\theta^0) + \Embb[\lambda^{\Top} r^0(\omega)]$ in expectation regardless of the retailer's net revenue target $F$.

The equivalent \emph{collective} welfare effects of DERs integrated under both models (characterized by Cor. \ref{cor:decentralized} and \ref{cor:centralized}) are in contrast to their \emph{individual} welfare effects.
As suggested above, welfare gains (or losses) from decentralized DERs are captured individually by DER-integrating customers as reductions in their bills and as bill reductions or increments for all other customers due to the adjustments to the connection charge $A^*_{\textsc{dec}}$.
This allocation of welfare gains constitutes an interclass cross-subsidy between customers.
Conversely, welfare gains from centralized DERs are uniformly captured by all customers as reductions in the connection charge $A^*_{\textsc{cen}}$.

Lastly, an implication of $T^*_{\textsc{cen}}$ is the likely impact it has on the adoption of DERs.
The reduction in the connection charge characterized by $A^*_{\textsc{cen}}$ relative to $A^*$ is the net benefit perceived by each customer due to the integrated centralized DERs.
Hence, customers should be willing to let the retailer integrate DERs even if they entail capital costs that offset a portion of said reductions in the connection charge.


\section{An Empirical Case Study} \label{sec:caseStudy2}

In this section, we analyze a case study of a hypothetical distribution utility that faces New York city's wholesale prices and residential demand for an average summer day.
We compare the performance of several day-ahead tariffs with hourly prices 
at different levels of solar and storage capacity.

Besides the optimal two-part tariff, we study other tariff structures with two pricing alternatives (flat pricing or hourly dynamic pricing) and with daily connection charges fixed at various levels: zero, a nominal value reflecting Con Edison's connection charge, the nominal value plus $0.33$ \$/day ($10$ \$/month), and the nominal value plus $1.66$ \$/day ($50$ \$/month).
Similar tariff reforms are being proposed in practice to solve utilities' fixed cost recovery problem \cite{LBNL16b}.
Given a tariff structure, we optimize the non-fixed parameters to maximize the expected consumer surplus subject to revenue adequacy.

This case study uses the same demand model as in Part I \cite{MunozTong16partIarxiv}, which comprises a linear demand function and a quadratic utility function for each customer.
We use publicly available energy sales data and rates from Con Edison for the 2015 Summer to fit the demand model.
Con Edison's default tariff for its $2.2$ million residential customers is essentially a two-part tariff $T^{\textsc{ce}}$ with a \emph{flat} price of $\pi^{\textsc{ce}}=17.2$ \cent/kWh and a connection charge of $15.76$ \$/month ($A^{\textsc{ce}}=0.53$ \$/day).
We use day-ahead wholesale prices for NYC from NYISO.

\subsection{Base case} \label{sec:caseStudy2:baseCase}

This is the case without DERs and nominal tariff $T^{\textsc{ce}}$.
Throughout the case study, we assume an average price elasticity of the total daily demand of $\overline{\varepsilon}(\pi^{\textsc{ce}})=-0.3$ at $\pi^{\textsc{ce}}$, which is a reasonable estimate of the short-term own-price elasticity of electricity demand \ifThesis{\cite{Borenstein05, Lijesen07, EPRI08}}{\cite{EPRI08}}.
Moreover, we consider a total of $M=2.2$ million residential customers and use the tariff $T^{\textsc{ce}}$ to compute the utility's average daily revenue from the residential segment and the portion that contributes towards fixed costs, which amount respectively to $\overline{\textsf{rev}}(T^{\textsc{ce}})= \$7.19$ million (M) dollars and 
$$
F^{\textsc{ce}} := \overline{\rp}_{0}(T^{\textsc{ce}}) = \overline{\textsf{rev}}(T^{\textsc{ce}}) - \Embb \left[ \lambda^{\Top} D(\b{1}\pi^{\textsc{ce}},\omega) \right]=\$5.83 \text{M.}
$$    
For the sake of brevity, the details of these computations already described in Part I \cite{MunozTong16partIarxiv} are not reproduced here.


\setlength\fboxsep{0pt}
\setlength\fboxrule{0pt}


We illustrate in Fig. \ref{fig:example:ParetoFrontBaseCase} the expected retailer surplus ($\overline{\rp}$) and expected consumer surplus ($\overline{\cs}$) induced by the revenue adequate tariff that maximizes the expected consumer surplus within each tariff structure for different values of $F$.
For each tariff structure, the resulting parametric curve is a Pareto front that quantifies the compromise between $\overline{\cs}$ and the $\overline{\rp}$ target, $F$.
We plot these curves as (possibly negative) surplus \emph{gains} relative to the values induced by $T^{\textsc{ce}}$, $\overline{\rp}_0(T^{\textsc{ce}})$ and $\overline{\cs}_0(T^{\textsc{ce}})=\$9.54$M, normalized by $\overline{\textsf{rev}}(T^{\textsc{ce}})$.

We make some observations from Fig. \ref{fig:example:ParetoFrontBaseCase}.
First, the $-1$ slope of the Pareto front associated to the optimal two-part tariff $T^*$ corroborates that the induced efficiency $\overline{\sw}(T^*)$ does not depend on $F$.
Conversely, the larger the $F$, the more inefficient the suboptimal two-part tariffs considered become.
This can be seen from the non-unitary slopes exhibited by all tariffs except $T^*$.
Second, at the nominal $\overline{\rp}$ target $F^{\textsc{ce}}$ (\ie the horizontal axis),
significant differences in the induced $\overline{\cs}$ gains are observed among the tariffs.
In particular, moving from flat prices to hourly prices improves $\overline{\cs}$ by approximately $1\%$ ($\$72$k/day).
A more significant $\overline{\cs}$ gain ($8.1\%$) is brought by also increasing the connection charge to the optimal level (which amounts to $A^{*}=2.65$ \$/day or $79.5$ \$/month).
Conversely, decreasing the connection charge to zero reduces $\overline{\cs}$ by $4.8\%$.
These empirical computations suggest that additional fixed costs can be recovered more efficiently by increasing connection charges than by pricing more dynamically.

\subsection{Tariff structure and net benefits of DERs} \label{sec:caseStudy2:efficiency}

We now analyze the combined impact of tariff structure and DER integration on consumers' surplus.
We measure changes in $\overline{\cs}$ relative to $\overline{\cs}_0(T^{\textsc{ce}})$ and normalized by $\overline{\textsf{rev}}(T^{\textsc{ce}})$.

\subsubsection{Customer-integrated DERs} \label{sec:caseStudy2:decentralized}

We start by estimating changes in $\overline{\cs}$ as a function of the solar and battery storage aggregate capacity integrated by customers.
The tariffs here considered are applied to the hourly net metered demand, so they differ from existing net metering tariffs with rolling credit.
Moreover, we model the integration of renewable resources using hourly solar PV generation data from a simulated 5kW-DC-capacity rooftop system located in NYC\footnote{``Typical year'' solar power data for the same months as temperature is taken from NREL's PVWatts Calculator available in \url{http://pvwatts.nrel.gov}.}.
Similarly, we consider the basic specifications of a $6.4$ kWh Tesla Powerwall battery\footnote{More precisely, a $3.3$ kW charging/discharging rate and a 96\% charging/discharging efficiency are used.}.
We integrate as many of these systems as necessary to reach the specified level of capacity.


\ifdefined\IEEEPARstart

\begin{figure*}[t]
\centering
\subcaptionbox{Pareto front for base case (zoom-in).\label{fig:example:ParetoFrontBaseCase}}
{\includegraphics[width=0.32\linewidth]{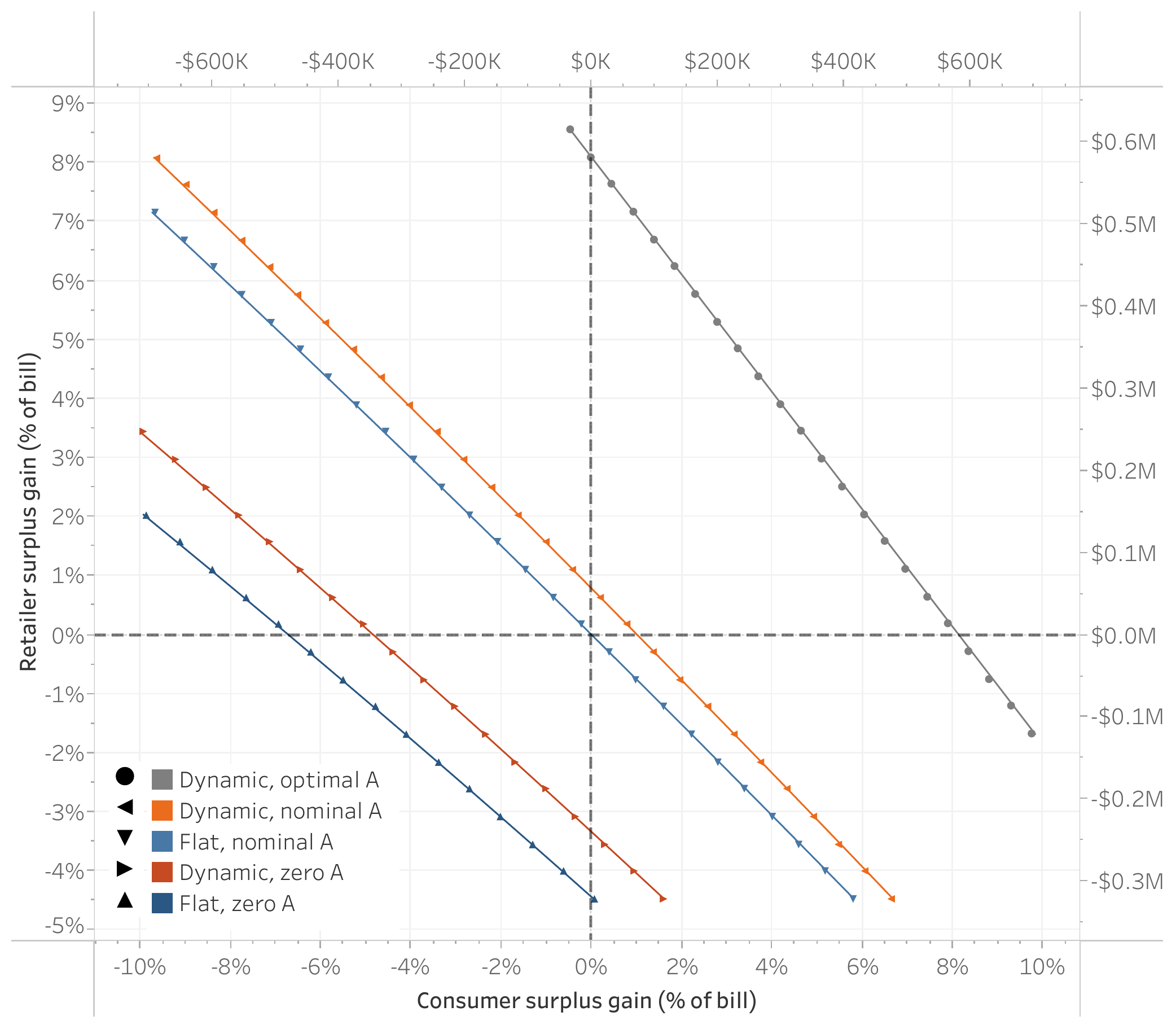}}
\subcaptionbox{Pareto fronts with decentralized DERs.\label{fig:example:PFdecentralized}}
{\includegraphics[width=0.32\linewidth]{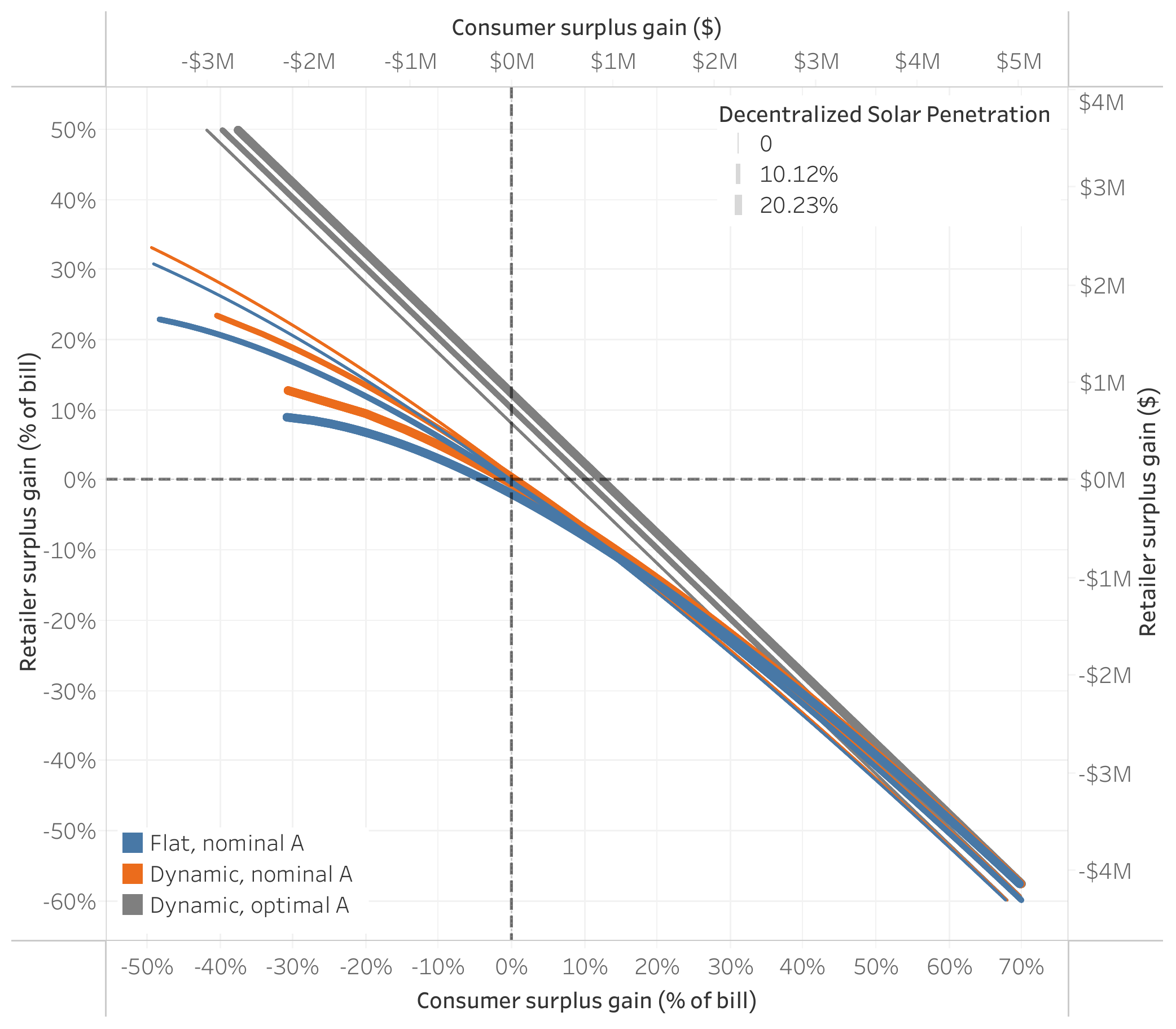}}
\subcaptionbox{Pareto fronts with centralized DERs.\label{fig:example:PFcentralized}}
{\includegraphics[width=0.32\linewidth]{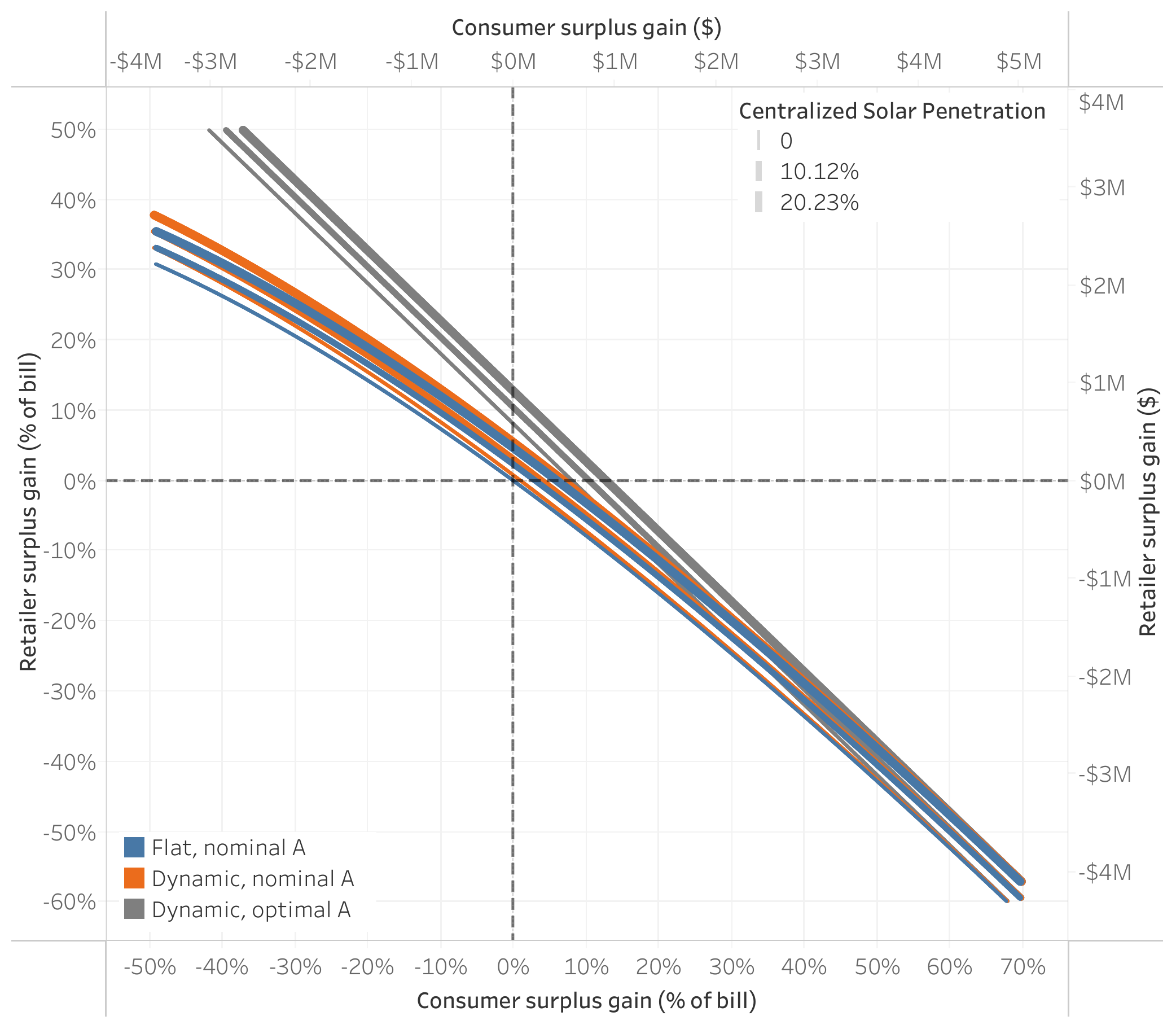}}
\caption{\protect\subref{fig:example:ParetoFrontBaseCase} Normalized retailer surplus target \textit{v.s.} induced consumer surplus gain (Pareto front) for various tariffs in base case (\ie no DERs). In \protect\subref{fig:example:PFdecentralized} and \protect\subref{fig:example:PFcentralized}, Pareto fronts for base case and two cases with different DER integration levels ($10.12\%$ and $20.23\%$) are compared.}
\label{fig:example:paretoFront}
\end{figure*}

	
\else

\begin{figure*}[t]
\centering
\subfloat[][Zoom into neighborhood of $\left(\overline{\cs}(T^{\textsc{ce}}),\overline{\rp}(T^{\textsc{ce}}) \right)$.]{
\includegraphics[width=1\linewidth]{ParetoFrontZoomIn}\label{fig:example:ParetoFrontBaseCase}}
\caption{Normalized retailer surplus target \textit{v.s.} induced consumer surplus gain (Pareto front) for various tariffs in base case (\ie no DERs).}
\label{fig:ParetoFronts:partII}
\end{figure*}

\begin{figure*}[t]
\centering
\subfloat[][Pareto fronts with decentralized DERs.]{\includegraphics[width=1\linewidth]{ParetoFrontDecentralized2}\label{fig:example:PFdecentralized}}\\
\subfloat[][Pareto fronts with centralized DERs.]{\includegraphics[width=1\linewidth]{ParetoFrontCentralized2}\label{fig:example:PFcentralized}}
\caption{\protect\subref{fig:example:ParetoFrontBaseCase} Normalized retailer surplus target \textit{v.s.} induced consumer surplus gain (Pareto front) for various tariffs in base case (\ie no DERs) and two cases with different DER integration levels ($10.12\%$ and $20.23\%$) are compared.}
\label{fig:example:paretoFront}
\end{figure*}

\fi

In Fig. \ref{fig:example:PFdecentralized} we plot the Pareto fronts associated to three types of tariffs
and three decentralized solar PV integration levels.
This figure is a zoomed-out version of Fig. \ref{fig:example:ParetoFrontBaseCase} computed for three DER integration levels.
As such, it gives a rough intuition of how decentralized DERs transform the Pareto fronts of different tariff structures and, in turn, affect $\overline{\cs}$.
In general, horizontal differences between the Pareto fronts represent changes in $\overline{\cs}$ due to tariff structure and/or to different levels of decentralized DER integration, for certain $F$.
Evidently, for any $F$, decentralized DERs bring $\overline{\cs}$ gains if the flat tariff structure is replaced by the optimal two-part tariff $T^*_{\textsc{dec}}$.
Conversely, $\overline{\cs}$ losses are brought by the DERs for $F=F^{\textsc{ce}}$ if the adjusted flat tariff structure is kept, or if it is replaced by the ``dynamic, nominal A'' structure.
We quantify these changes in $\overline{\cs}$ for $F=F^{\textsc{ce}}$ explicitly with the following parametric analysis over the decentralized DER integration level.

In Fig. \ref{fig:example:CSgains:a}, for several tariff structures, we plot the normalized $\overline{\cs}$ gains (or, equivalently, the $\overline{\sw}$ gains) caused by increments in the PV capacity integrated by customers and the corresponding updates to the tariff required to maintain revenue adequacy, \ie $\overline{\rp}(T)=F^{\textsc{ce}}$.
This case assumes that the storage capacity integrated is half the PV capacity.

In particular, Fig. \ref{fig:example:CSgains:a} shows how integrating decentralized DERs can trigger both efficiency gains and losses depending on the tariff structure.
For example, the curve for the adjusted flat tariff $T^{\textsc{ce}}$ suggests that maintaining revenue adequacy with flat rate increments would cause DERs to bring no significant net gains or losses in $\overline{\cs}$ and $\overline{\sw}$ for small levels of integration.
However, DER integration levels beyond $500$ MW would bring increasingly larger losses in $\overline{\cs}$ and $\overline{\sw}$ of $1.3\%$ at $1.1$ GW and $15\%$ at $2$ GW.
A similar performance is shown by the optimal linear tariff (dynamic pricing) with nominal connection charge $A^{\textsc{ce}}$.
The gain of $1\%$ it exhibits with no DERs vanishes to $0\%$ at $1.1$ GW of PV, becoming net efficiency losses for higher levels of DER.
The optimal linear tariffs with higher connection charges ---$10$ and $50$ \$/month larger than the nominal---, which bring initial efficiency gains of $2.8\%$ and $7.7\%$, respectively.
While the gains of the former increase to then decrease after reaching a maximum of $3.5\%$, the gains of the latter monotonically increase reaching $13.6\%$ at the maximum PV capacity considered, $2.2$ GW.
Other example is the flat tariff with no connection charge, which starts with efficiency losses of $6.8\%$ that increase sharply up to $20.4\%$ at $1.1$ GW.
Lastly, the optimal two-part tariff starts with an efficiency gain of $8.2\%$, and it lets customer-integrated DERs to generate their full value, which amounts to an additional $6.6\%$ of efficiency at $2.2$ GW or $3\%$ per GW.
In other terms, the efficiency gains foregone by using the adjusted flat tariff $T^{\textsc{ce}}$ rather than the optimal two-part tariff $T^*$ increase linearly with the level of DER integration and reach $29.3\%$ (or $\$2.11$ M/day) at $2$ GW. 

In summary, connection charges embody a method for fixed cost recovery that seems to be even more effective than dynamic pricing in the sense that it can harness at least $90\%$ of the efficiency gains attained by the optimal two-part tariff for all the integration levels considered, whereas dynamic pricing alone can harnesses at most $12.5\%$ and it generates efficiency losses for higher integration levels.
A word of caution on tariffs with high connection charges and lower flat prices, however, is that they induce customers to consume more on peak, thus precipitating the need for network upgrades that increase the retailer's fixed costs in a way not captured by our model.
This problem can be tackled using dynamic prices and high connection charges to recover fixed costs, such as the optimal two-part tariff, and more forcefully by considering endogenous fixed costs dependent of the coincident peak net-load. 

\setlength\fboxsep{0pt}
\setlength\fboxrule{0pt}

\ifdefined\IEEEPARstart

\begin{figure*}[t]
\centering
\includegraphics[width=0.99\linewidth]{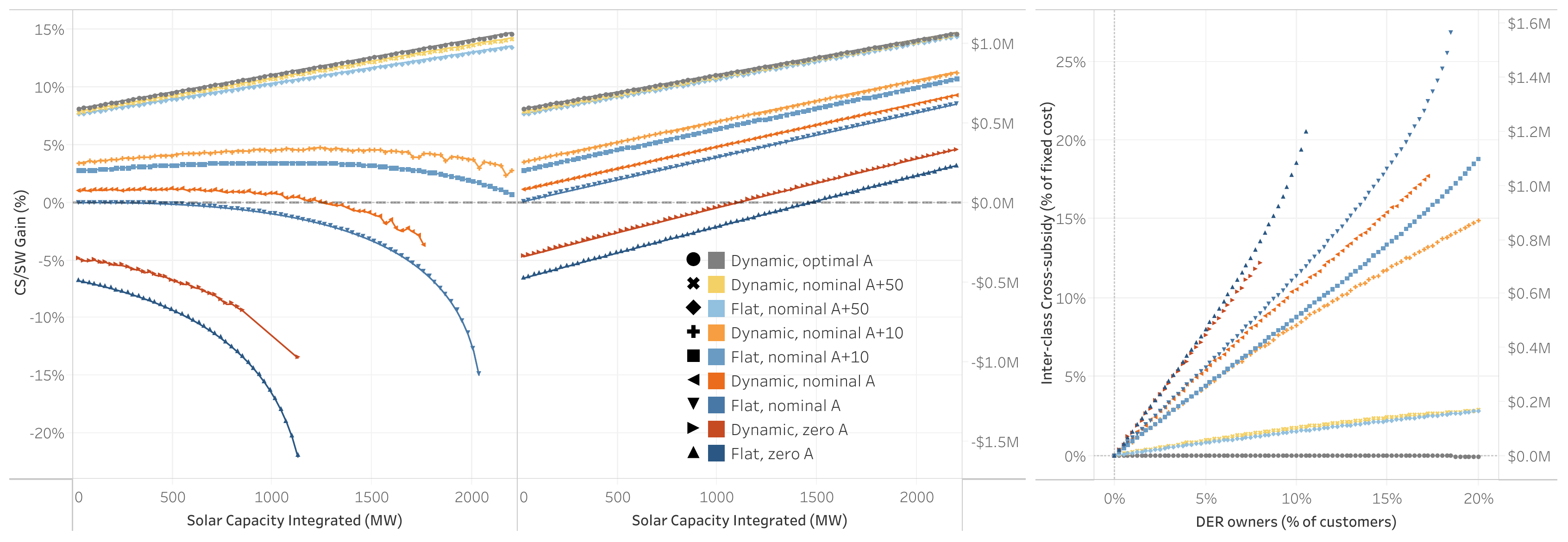}\\[-15pt]
\subcaptionbox{Decentralized DER integration.\label{fig:example:CSgains:a}}
{\rule{0.3\linewidth}{0pt}}
\subcaptionbox{Centralized DER integration.\label{fig:example:CSgains:b}}
{\rule{0.3\linewidth}{0pt}}
\subcaptionbox{Cross-subsidies.\label{fig:example:cross-subsidies}}
{\rule{0.3\linewidth}{0pt}}
\caption{Expected gains in consumer and social surplus induced by \protect\subref{fig:example:CSgains:a} behind-the-meter solar-plus-battery capacity and \protect\subref{fig:example:CSgains:b} retailer-integrated solar-plus-battery capacity under different types of tariffs.
Gains are measured relative to base case with tariff $T^{\textsc{ce}}$ and no DERs.
In \protect\subref{fig:example:cross-subsidies}, cross-subsidies from customers without solar to customers with solar \textit{v.s.} level of behind-the-meter solar integration.}
\label{fig:example:CSgains}
\end{figure*}

	
\else

\begin{figure}[h]
\centering
\subfloat[][Decentralized DER integration.]{
\includegraphics[width=0.99\linewidth]{CSoverallGainsDec2}\label{fig:example:CSgains:a}}\\
\subfloat[][Centralized DER integration.]{
\includegraphics[width=0.99\linewidth]{CSoverallGainsCen}\label{fig:example:CSgains:b}}
\caption{Expected gains in consumer and social surplus induced by \protect\subref{fig:example:CSgains:a} behind-the-meter solar-plus-battery capacity and \protect\subref{fig:example:CSgains:b} retailer-integrated solar-plus-battery capacity under different types of tariffs.
Gains are measured relative to base case with tariff $T^{\textsc{ce}}$ and no DERs.}
\label{fig:example:CSgains}
\end{figure}

\begin{figure}[h]
\centering
\includegraphics[width=0.7\linewidth]{crossSubsidiesTriangular}
\caption{Cross-subsidy from customers without solar to customers with solar \textit{v.s.} level of behind-the-meter solar integration.}
\label{fig:example:cross-subsidies}
\end{figure}

\fi

%
%

\subsubsection{Retailer-integrated DERs} \label{sec:caseStudy2:centralized}

We now estimate changes in surplus as a function of the solar and battery storage capacity integrated by the retailer.
For the sake of a fair comparison, DERs with the same characteristics as before are used.

In Fig. \ref{fig:example:CSgains:b}, for several tariff structures, we plot the normalized gains in $\overline{\cs}$ and $\overline{\sw}$ caused by increments in the PV and storage capacity integrated (in a 2:1 ratio) by the retailer and the corresponding tariff updates required to maintain revenue adequacy.
The figure reveals that the DERs bring monotonic surplus gains under all the tariffs considered.
This is because the benefits brought by centralized DERs are not offset by the consumption inefficiencies induced by decentralized DERs.
In fact, the inefficiencies induced by suboptimal tariffs without DERs are slightly mitigated by centralized DERs because these DERs help the retailer recover a small portion of $F$.

Hence, Fig. \ref{fig:example:CSgains:b} suggests that the changes in $\overline{\cs}$ and $\overline{\sw}$ brought by retailer-integrated DERs are virtually unbiased by tariff structure.
This is unlike customer-integrated DERs whose effect on $\overline{\cs}$ and $\overline{\sw}$ is significantly biased by tariff structure changes, and specially, by the reliance on retail price markups for fixed cost recovery, as it is evident in Fig. \ref{fig:example:CSgains:a}.
In other words, under tariffs with significant retail markups, while centralized DERs generally bring surplus gains, decentralized DERs tend to mitigate surplus gains or bring surplus losses.
It also clear from Fig. \ref{fig:example:CSgains} that dynamic pricing and higher connection charges help consistently (\ie regardless the level of DER integration) to mitigate existing inefficiencies.
Notably, connection charges seem to offer a much more effective measure to mitigate such inefficiencies than dynamic pricing.


\subsection{Cross-subsidies induced by net metering} \label{sec:caseStudy2:xsubsidies}

Considering the inequity concerns raised by using net metering as a mechanism to compensate DERs \ifThesis{\cite{EidEtal14,BorlickWood14}\cite[Sec. 9.5]{MIT15}}{\cite{EidEtal14}\cite[Sec. 9.5]{MIT15}}, it is instructive to quantify the cross-subsidies induced by the tariffs in the previous section.
To that end, we compute the cross-subsidies between PV owners and non-PV owners for different levels of customer-integrated solar PV capacity.

For a given tariff structure and level of (decentralized) solar integration, the cross-subsidy is computed by first obtaining the contribution that PV owners make towards the fixed costs $F$.
Similarly, the contribution that PV owners would make under a version of the given tariff that settles consumption and generation separately is also computed.
This version, which is also optimized subject to revenue adequacy, settles all consumption at the rates $\pi$ and all generation at the prices $\overline{\lambda}$\footnote{Customer-integrated storage is not considered in this analysis because nonlinear tariffs make customers' problem fundamentally more complicated.
}.
Clearly, PV owners contribute less towards $F$ under the net metering tariff than its counterpart if the associated prices are marked up as a means to recover the fixed costs $F$.
We compute said cross-subsidy as the difference between these two values, normalized by $F$.
Intuitively, cross-subsidies are the difference between the costs that each group \emph{should} pay for and those they \emph{actually} pay for, due to net metering.

The computation of cross-subsidies requires specifying individual demand functions and  the distribution of solar capacity between customers.
We consider a simple illustrative case.
Customers have identical demand functions except for a scaling parameter $\sigma_i$ satisfying $\sigma_i=i \cdot \sigma$ for some $\sigma > 0$.
The 5-kW solar installations are allocated to the largest consumers, who have the greatest incentive to invest in solar generation.
The resulting inter-class cross-subsidies are depicted in Fig. \ref{fig:example:cross-subsidies}.
%
Evidently, all tariffs induce non-trivial cross-subsidies except for the optimal two-part tariff which yields virtually no cross-subsidy (in spite of the discussion in Sec. \ref{sec:performance with DER}).
This is not entirely surprising since pricing according to $\pi^*_{\textsc{dec}}=\overline{\lambda}$ is efficient and consistent with cost causality.
Hence, cross-subsidies with such tariff happen only through the second term in \eqref{eq2:A decentralized} which is rather small compared to $A^*$.
The cross-subsidies of all the other tariffs increase with PV capacity, and they do it at an increasingly faster rate for flat tariffs.

\section{Conclusions} \label{sec:conclusions2}


We leverage the analytical framework developed in \cite{MunozTong16partIarxiv} to study how retail electricity tariff structure can distort the net benefits brought by DERs integrated by customers and their retailer.
This work is an application of Ramsey pricing with extensions to accommodate the integration of DERs.


Our analysis offers several conclusions.
First, while net metering tariffs that rely on flat and higher prices to maintain revenue adequacy provide increasingly stronger incentives for customers to integrate renewables%
, they induce increasingly larger cross-subsidies and consumption inefficiencies that can outweigh renewables' benefits.
These significant inefficiencies have draw little attention in the literature compared to the cross-subsidies.
Second, net metering tariffs can achieve revenue adequacy without compromising efficiency by using marginal-cost-based dynamic prices and higher connection charges.
These tariffs, however, provide little incentive to integrate renewables.
Third, retailer-integrated DERs bring customers net benefits that are less dependent on tariff structure, and they cause no tariff feedback loops.
As such, this alternative to behind-the-meter DERs seems worth exploring.

This study represents an initial point of analysis, for it has various limitations.
First, policy objectives beyond efficiency and revenue sufficiency ---often considered in practice--- are here ignored.
Practical criteria such as bill stability make ``desirable'' tariffs hard to be ever attained.
Second, customer disconnections are assumed to be not plausible as a customer choice.
This assumption becomes increasingly less realistic with the decline of DER costs.
Lastly, retailer non-energy costs are assumed to be fixed and independent of the coincident (net) peak load.
Relaxing this assumption leads to peak-load pricing formulations \cite{CrewEtAl95}.
We discuss one such relaxation in \cite{Munoz16}, where capacity costs are recovered with a \emph{demand charge} applied to the net demand coincident with the peak period.




\bibliographystyle{IEEEtran}

\begin{thebibliography}{10}
\providecommand{\url}[1]{#1}
\csname url@samestyle\endcsname
\providecommand{\newblock}{\relax}
\providecommand{\bibinfo}[2]{#2}
\providecommand{\BIBentrySTDinterwordspacing}{\spaceskip=0pt\relax}
\providecommand{\BIBentryALTinterwordstretchfactor}{4}
\providecommand{\BIBentryALTinterwordspacing}{\spaceskip=\fontdimen2\font plus
\BIBentryALTinterwordstretchfactor\fontdimen3\font minus
  \fontdimen4\font\relax}
\providecommand{\BIBforeignlanguage}[2]{{%
\expandafter\ifx\csname l@#1\endcsname\relax
\typeout{** WARNING: IEEEtran.bst: No hyphenation pattern has been}%
\typeout{** loaded for the language `#1'. Using the pattern for}%
\typeout{** the default language instead.}%
\else
\language=\csname l@#1\endcsname
\fi
#2}}
\providecommand{\BIBdecl}{\relax}
\BIBdecl

\bibitem{MunozTong16partIarxiv}
D.~Munoz-Alvarez and L.~Tong, ``{On the Efficiency of Dynamic Retail Tariffs
  with Connection Charges---Part I: A Stochastic Optimization Formulation},''
  \emph{ArXiv Preprint}, 2017.

\bibitem{BrownSibley86}
S.~J. Brown and D.~S. Sibley, \emph{The theory of public utility
  pricing}.\hskip 1em plus 0.5em minus 0.4em\relax Cambridge University Press,
  1986.

\bibitem{MunozTong16partIIarxiv}
D.~Munoz-Alvarez and L.~Tong, ``{On the Efficiency of Dynamic Retail Tariffs
  with Connection Charges---Part II: Distributed Renewable and Storage
  Resources},'' \emph{ArXiv Preprint}, 2017.

\bibitem{Munoz16}
D.~{Mu\~noz-\'Alvarez}, ``Regulatory approaches to the integration of renewable
  and storage resources into electricity markets,'' Ph.D. dissertation, Cornell
  University, 2017.

\bibitem{RenesesRodriguez14}
J.~Reneses and M.~P. Rodr\'iguez-Ortega, ``{Distribution pricing: theoretical
  principles and practical approaches},'' \emph{IET Generation, Transmission
  {\&} Distribution}, no.~8, pp. 1645--1655, oct 2014.

\bibitem{LBNL16b}
L.~Wood, J.~Howat, R.~Cavanagh, and S.~Borenstein, ``Recovery of utility fixed
  costs: Utility, consumer, environmental and economist perspectives,'' LBNL,
  Tech. Rep., 2016.

\bibitem{Costello15}
K.~W. Costello, ``{Major challenges of distributed generation for state utility
  regulators},'' \emph{Electricity Journal}, vol.~28, no.~3, pp. 8--25, 2015.

\bibitem{NREL13}
L.~Bird, J.~Mclaren, and J.~Heeter, ``{Regulatory Considerations Associated
  with the Expanded Adoption of Distributed Solar},'' NREL, Tech. Rep.
  November, 2013.

\bibitem{MIT15}
R.~Schmalensee and V.~Bulovic, ``{The Future of Solar Energy},'' {Massachusetts
  Institute of Technology}, Tech. Rep., 2015.

\bibitem{NARUC16}
{National Association of Regulatory Utility Commissioners (NARUC)},
  ``{Distributed Energy Resources Rate Design and Compensation},'' {National
  Association of Regulatory Utility Commissioners}, Manual November, 2016.

\bibitem{EidEtal14}
C.~Eid, J.~{Reneses Guill{\'{e}}n}, P.~{Fr{\'{i}}as Mar{\'{i}}n}, and
  R.~Hakvoort, ``{The economic effect of electricity net-metering with solar
  PV: Consequences for network cost recovery, cross subsidies and policy
  objectives},'' \emph{Energy Policy}, vol.~75, pp. 244--254, 2014.

\bibitem{JargstorfBelmans15}
J.~Jargstorf and R.~Belmans, ``{Multi-objective low voltage grid tariff
  setting},'' \emph{IET Generation, Transmission {\&} Distribution}, vol.~9,
  no.~15, pp. 2328--2336, 2015.

\bibitem{DarghouthEtal16}
N.~R. Darghouth, R.~H. Wiser, G.~Barbose, and A.~D. Mills, ``{Net metering and
  market feedback loops: Exploring the impact of retail rate design on
  distributed PV deployment},'' \emph{Applied Energy}, vol. 162, 2016.

\bibitem{Sioshansi16}
R.~Sioshansi, ``{Retail electricity tariff and mechanism design to incentivize
  distributed renewable generation},'' \emph{Energy Policy}, vol.~95, 2016.

\bibitem{ChewEtal12}
M.~Chew, M.~Heling, C.~Kerrigan, D.~Jin, A.~Tinker, M.~Kolb, S.~Buller, and
  L.~Huang, ``Modeling dg adoption using electric rate feedback loops,'' in
  \emph{31st USAEE/IAEE North American Conference}, 2012.

\bibitem{CaiEtal13}
D.~W. Cai, S.~Adlakha, S.~H. Low, P.~De~Martini, and K.~M. Chandy, ``Impact of
  residential pv adoption on retail electricity rates,'' \emph{Energy Policy},
  vol.~62, pp. 830--843, 2013.

\bibitem{Kind13}
P.~Kind, ``{Disruptive Challenges: Financial Implications and Strategic
  Responses to a Changing Retail Electric Business},'' Edison Electric
  Institute, Tech. Rep. January, 2013.

\bibitem{MIT13}
A.~Bharatkumar, S.~Burger, J.~D. Jenkins, J.~L. Dantec, I.~J.
  P{\'{e}}rez-Arriaga, R.~D. Tabors, and C.~Batlle, ``{The MIT Utility of the
  Future Study},'' MIT, Tech. Rep. December, 2013.

\bibitem{CostelloHemphill14}
K.~W. Costello and R.~C. Hemphill, ``{Electric utilities' 'death spiral':
  Hyperbole or reality?}'' \emph{Electricity Journal}, vol.~27, no.~10, 2014.

\bibitem{RMI14}
\BIBentryALTinterwordspacing
{Rocky Mountain Institute}, ``The economics of grid defection,'' Tech. Rep.,
  2014. [Online]. Available: \url{http://www.rmi.org/}
\BIBentrySTDinterwordspacing

\bibitem{RMI15}
\BIBentryALTinterwordspacing
------, ``The economics of load defection,'' Tech. Rep., 2015. [Online].
  Available: \url{http://www.rmi.org/}
\BIBentrySTDinterwordspacing

\bibitem{ChenEtal12}
L.~Chen, N.~Li, L.~Jiang, and S.~H. Low, ``Optimal demand response: problem
  formulation and deterministic case,'' in \emph{Control and optimization
  methods for electric smart grids}.\hskip 1em plus 0.5em minus 0.4em\relax
  Springer, 2012, ch.~3, pp. 63--85.

\bibitem{TangEtal14}
W.~Tang, R.~Jain, and R.~Rajagopal, ``Stochastic dynamic pricing: Utilizing
  demand response in an adaptive manner,'' in \emph{Decision and Control (CDC),
  IEEE Conference on}, 2014, pp. 6446--6451.

\bibitem{JiaTong16b}
L.~Jia and L.~Tong, ``Renewables and storage in distribution systems:
  Centralized vs. decentralized integration,'' \emph{IEEE J. on Sel. Areas in
  Comm.}, vol.~34, no.~3, pp. 665--674, March 2016.

\bibitem{HanifEtal16}
S.~Hanif, T.~Massier, T.~Hamacher, and T.~Reindl, ``{Evaluating demand response
  in the presence of solar PV: Distribution grid perspective},'' in \emph{Smart
  Energy Grid Engineering, IEEE}.\hskip 1em plus 0.5em minus 0.4em\relax IEEE,
  aug 2016, pp. 392--397.

\bibitem{JiaTong16a}
L.~Jia and L.~Tong, ``Dynamic pricing and distributed energy management for
  demand response,'' \emph{{Smart Grid, IEEE Trans. on}}, vol.~7, 2016, earlier
  arXiv version: \url{http://arxiv.org/abs/1601.02319}.

\bibitem{EPRI08}
{Electric Power Research Institute (EPRI)}, ``{Price Elasticity of Demand for
  Electricity : A Primer and Synthesis},'' Tech. Rep., 2008.

\bibitem{CrewEtAl95}
M.~A. Crew, C.~S. Fernando, and P.~R. Kleindorfer, ``{The theory of peak-load
  pricing: A survey},'' \emph{Journal of Regulatory Economics}, vol.~8, no.~3,
  pp. 215--248, nov 1995.

\end{thebibliography}

\appendix
\label{sec:appendix:part2}

\ifdefined\IEEEPARstart

\begin{figure*}
\begin{align}
\overline{\cs}^i(T) &= \max_{q^i(\cdot),s^i(\cdot)} \bigg\{\Embb \left[ S^i(q^i(\omega^i),\omega^i) - T\Big(q^i(\omega^i)-r^i(\omega^i)-s^i(\omega^i)\Big) \right] \Big| s^i(\omega^i) \in \Umsc(\theta^i) \bigg\} \nn \\
&= \Embb \bigg[ S^i(q^{i*}(T,\omega^i),\omega^i) - T\Big( \underbrace{q^{i*}(T,\omega^i)-r^{i}(\omega^i)-s^{i*}(T,\omega^i)}_{d^{i*}(T,\omega^i)} \Big) \bigg]
\label{eq2:cs with DERs long}
\end{align}
\hrulefill
\end{figure*}

\else
\fi

\subsubsection*{Proof of Theorem \ref{thm:decentralized}} \label{proof:thm:decentralized}

To solve the optimization problem \eqref{eq2:reg problem} over affine tariffs of the form $T(d)=A+\pi^{\Top}d$, we obtain expressions for $\overline{\cs}(T)$ and $\overline{\rp}(T)$ in terms of the parameters $\pi$ and $A$, considering the customer-integrated DERs.

On one hand, from the separability implications of the linearity of $T$ on the customers' problem established in Sec. \ref{sec:decentralizedIntegration}, we have that customers with DERs obtain an expected surplus
\begin{align*}
\overline{\cs}^i(T) &= \Embb[S^i(D^i(\pi,\omega^i))-T(d^i(\pi,\omega^i))] \\
					&= \overline{\cs}^i_0(T) + \Embb[\pi^{\Top} (s^*(\pi,\theta^i) + r^i(\omega^i))],
\end{align*}
where $\overline{\cs}^i_0(T)$, the expected consumer surplus without DERs, is computed according to \eqref{eq2:cs}, \ie
\begin{align} \label{eq2:cs without DERs}
\overline{\cs}^i_0(T) = \Embb \big[ S^i(D^i(\pi,\omega^i),\omega^i) - \pi^{\Top}D^i(\pi,\omega^i) \big] - A.
\end{align}
On the other hand, since this case does not consider retailer-integrated DERs, the retailer derives an expected surplus
\begin{align*}
\overline{\rp}(T)	&= \sum_{i=1}^M \Embb[ T(d^i(\pi,\omega^i)) - \lambda^{\Top}d^i(\pi,\omega^i) ] \\
					&= \overline{\rp}_0(T) - \sum_{i=1}^M \Embb[ (\pi-\lambda)^{\Top} ( s^*(\pi,\theta^i) + r^i(\omega^i) ) ]
\end{align*}
where $\overline{\rp}_0(T)$, the expected retailer surplus without DERs, is computed according to \eqref{eq2:rs}, \ie
$$
\overline{\rp}_0(T) = A\cdot M + \sum_{i=1}^M \Embb \left[ (\pi - \lambda)^{\Top}D^i(\pi,\omega^i) \right].
$$

We can now solve the constraint $\overline{\rp}(T)=F$ for $A$,
\begin{align} \label{eq2:connection charge decentralized}
A=\mbox{$\frac{1}{M}$} \left( F - \Embb \left[ \mbox{$\sum_{i=1}^M$} (\pi-\lambda)^{\Top} d^i(\pi,\omega^i) \right] \right)
\end{align}
and replace it in the objective function of problem \eqref{eq2:reg problem}, $\overline{\cs}(T)$, which yields
\begin{align} \label{eq2:cs decentralized}
\overline{\cs}(T) 	&= \overline{\sw}_0(T) - F + \sum_{i=1}^M \Embb[ \lambda^{\Top} ( s^*(\pi,\theta^i) + r^i(\omega^i) ) ] 
\end{align}
where $\overline{\sw}_0(T) = \overline{\cs}_0(T) + \overline{\rp}_0(T)$, the expected total surplus that $T$ would induce in the absence of DERs, is given by
\begin{align} \label{eq2:sw no DERs}
\overline{\sw}_0(T)	&= \sum_{i=1}^{M} \Embb\left[ S^i(D^{i}(\pi,\omega^i), \omega^i) -  \lambda^{\Top} D^{i}(\pi,\omega^i) \right].
\end{align}

One can then show that $\overline{\sw}_0(T)$ is maximized over $\pi$ at $\pi_{\textsc{dec}} = \overline{\lambda}$ (see proof of Theorem 1 in \cite{MunozTong16partIarxiv}), when $\nabla_{\pi}D(\pi,\omega)$ and $\lambda$ are uncorrelated (see Cor. 1 in \cite{MunozTong16partIarxiv}).
Moreover, each term $\Embb[ \lambda^{\Top} s^*(\pi,\theta^i)]$ is also maximized over $\pi$ at $\pi_{\textsc{dec}} = \overline{\lambda}$ since
$$
\Embb[ \lambda^{\Top} s^*(\pi,\theta^i)] = \overline{\lambda}^{\Top} s^*(\pi,\theta^i) \leq \overline{\lambda}^{\Top} s^*(\overline{\lambda},\theta^i) = V(\overline{\lambda},\theta^i)
$$
for all $\pi \in \Rmbb^N$.
Consequently, the expression $\overline{\cs}(T)$ in \eqref{eq2:cs decentralized} is maximized over $\pi$ at $\pi^*_{\textsc{dec}} = \overline{\lambda}$.
The strict concavity of $\overline{\sw}_0(T)$ in $\pi$ (Prop. 4 in \cite{MunozTong16partIarxiv}) guarantees the uniqueness and optimality of $\pi^*_{\textsc{dec}}$ and $A^*_{\textsc{dec}}$.
Replacing $\pi$ and $d^i(\pi,\omega^i)$ in \eqref{eq2:connection charge decentralized} for $\overline{\lambda}$ and $d^i(\overline{\lambda},\omega^i)$ according to \eqref{eq2:demand with DERs}, respectively, yields the expression for $A^*_{\textsc{dec}}$ in \eqref{eq2:A decentralized}, where
$$
A^* = \mbox{$\frac{1}{M}$} \left( F + \tr \left( \cov \left( \lambda,D \left( \overline{\lambda},\omega \right) \right) \right) \right) 
$$
is the optimal connection charge without DERs for uncorrelated $\nabla_{\pi}D(\pi,\omega)$ and $\lambda$ (Corollary 1 in \cite{MunozTong16partIarxiv}).
\hfill $\blacksquare$

\subsubsection*{Proof of Corollary \ref{cor:decentralized}} \label{proof:cor:decentralized}

Theorem \ref{thm:decentralized} implies that if $\nabla_{\pi}D(\pi,\omega)$ and $\lambda$ are uncorrelated then $\pi^*_{\textsc{dec}}=\overline{\lambda}$, which does not depend on the parameter $F$.
Hence, it is clear from \eqref{eq2:sw no DERs} that
$
\overline{\sw}^*_0 \equiv \overline{\sw}_0(T^*),
$
where $T^*(q)=A^*+\pi^{*\Top}q$, 
is also independent from the parameter $F$.
It follows from \eqref{eq2:cs decentralized} that
$$
\overline{\cs}(T^*_{\textsc{dec}})= \overline{\sw}^*_0 - F + \mbox{$\sum_{i=1}^{M}$} V(\overline{\lambda},\theta^i) + \Embb[ \lambda^{\Top} r^i(\omega^i)]
$$
and, since $\overline{\rp}(T^*_{\textsc{dec}})=F$ at optimality,
\begin{align*}
\overline{\sw}(T^*_{\textsc{dec}}) = \overline{\sw}^*_0 + \mbox{$\sum_{i=1}^{M}$} V(\overline{\lambda},\theta^i) + \Embb[ \lambda^{\Top} r^i(\omega^i)].
\end{align*}
Thus, $\overline{\sw}(T^*_{\textsc{dec}})$ is also independent of the parameter $F$.
\hfill $\blacksquare$

\subsubsection*{Proof of Theorem \ref{thm:optimalityDERs}} \label{proof:thm:optimalityDERs}

We show this result using arguments analogous to those in the proof of Theorem 3 in \cite{MunozTong16partIarxiv}.
That is, we show that the optimal two-part tariff $T^*_{\textsc{dec}}$ attains an upper bound for the performance of all \emph{ex-ante} tariffs derived from the social planner's problem.
To obtain a tight upper bound for ex-ante tariffs only (rather than a looser bound for all possibly ex-post tariffs), the social planner makes customers' decisions relying only on the information observable by each of them (\ie $\omega^i$) as opposed to based on global information (\eg $\xi=(\lambda,\omega^1,\ldots,\omega^M)$).

Consider the social planner's problem
\begin{subequations} \label{eq2:social planner with DERs}
\begin{align}
\max_{\{q^i(\cdot),s^i(\cdot)\}_{i=1}^M} &\ \overline{\sw} \\
							\text{s.t} \ \ &\ s^i(\omega^i) \in \Umsc(\theta^i), \qquad i=1,\ldots,M.
\end{align}
\end{subequations}
with
$$
\overline{\sw} = \Embb_{\xi} \left[ \sum_{i=1}^M S^i(q^i(\omega^i),\omega^i) - \lambda^{\Top}d^i(\omega^i) \right]
$$
which is related to the regulator's problem \eqref{eq2:reg problem} with customer-integrated DERs.
The notations $q^i(\omega^i)$, $s^i(\omega^i)$, and $d^i(\omega^i)$ indicate the restriction of the social planner to make (causal) decisions \emph{contingent only} on the local state of each customer $\omega^i$.
Recall from \eqref{eq2:cs with DERs} that the expected consumer surplus for a given ex-ante tariff is given by $\overline{\cs}(T) = \sum_{i=1}^M \overline{\cs}^i(T)$, where $\overline{\cs}^i(T)$ can be written as
\ifdefined\IEEEPARstart
 in \eqref{eq2:cs with DERs long},
\else
\begin{align}
\overline{\cs}^i(T) &= \max_{q^i(\cdot),s^i(\cdot)} \bigg\{\Embb \left[ S^i(q^i(\omega^i),\omega^i) - T\Big(q^i(\omega^i)-r^i(\omega^i)-s^i(\omega^i)\Big) \right] \Big| s^i(\omega^i) \in \Umsc(\theta^i) \bigg\} \nn \\
&= \Embb \bigg[ S^i(q^{i*}(T,\omega^i),\omega^i) - T\Big( \underbrace{q^{i*}(T,\omega^i)-r^{i}(\omega^i)-s^{i*}(T,\omega^i)}_{d^{i*}(T,\omega^i)} \Big) \bigg],
\label{eq2:cs with DERs long}
\end{align}
\fi
and note that the corresponding expected retailer surplus is given by
\begin{align*}
\overline{\rp}(T) = \Embb \left[ \sum_{i=1}^M T(d^{i*}(T,\omega^i)) - \lambda^{\Top} d^{i*}(T,\omega^i) \right],
\end{align*}
and the expected total surplus by 
\begin{align*}
\overline{\sw}(T) &= \Embb \bigg[ \sum_{i=1}^M S^i(q^{i*}(T,\omega^i),\omega^i) - \lambda^{\Top} d^{i*}(T,\omega^i) \bigg] \\
&= \Embb \bigg[ \sum_{i=1}^M S^i(q^{i*}(T,\omega^i),\omega^i) - \lambda^{\Top} q^{i*}(T,\omega^i)\bigg]  \\
& \qquad\qquad + \Embb \bigg[ \sum_{i=1}^M \lambda^{\Top} (r^i(\omega^i) + s^{i*}(T,\omega^i))\bigg].
\end{align*}

The following sequence of equalities/inequalities shows that problem \eqref{eq2:social planner with DERs} provides an upper bound to problem \eqref{eq2:reg problem}.
\begin{align} 
& \max_{T(\cdot)} \{ \overline{\cs}(T) \ | \ \overline{\rp}(T) = F \} + F \nn\\ 
&= \max_{T(\cdot)} \{ \overline{\cs}(T) + \overline{\rp}(T) \ | \ \overline{\rp}(T) = F \} \nn\\ 
&= \max_{T(\cdot)} \{ \overline{\sw}(T) \ | \ \overline{\rp}(T) = F \} \nn\\ 
&\leq \max_{T(\cdot)} \ \overline{\sw}(T) \label{eq2:reg problem UB 1} \\
&\leq \max_{\{q^i(\cdot)\}_{i=1}^M} \Embb\left[ \sum_{i=1}^MS^i(q^i(\omega^i),\omega^i) - \lambda^{\Top}q^i(\omega^i) \right] \nn\\
 &~~~~ +
\sum_{i=1}^M \max_{s^i(\cdot)} \Big\{ \left. \Embb\left[ \lambda^{\Top}(r^i(\omega^i) + s^i(\omega^i)) \right] \ \right| \  s^i(\omega^i) \in \Umsc(\theta^i) \Big\} \label{eq2:reg problem UB 2}  \\
&= \max_{\{q^i(\cdot),s^i(\cdot)\}_{i=1}^{M}} \ \Big\{ \ \overline{\sw} \ \big| \ s^i(\omega^i) \in \Umsc(\theta^i), \ i=1,\ldots,M \Big\}. \label{eq2:reg problem UB 3}
\end{align}

In particular, the inequality in \eqref{eq2:reg problem UB 2} holds because $\overline{\sw}(T)$ depends on $T$ only through $q^{i*}(T,\omega^i)$ and $s^{i*}(T,\omega^i)$.
This implies that maximizing $\overline{\sw}(T)$ directly over $\{q^i(\cdot),s^i(\cdot)\}_{i=1}^M$ rather than indirectly over $T(\cdot)$ is a relaxation of the optimization in \eqref{eq2:reg problem UB 1}.
Clearly, the problem in \eqref{eq2:reg problem UB 2} corresponds to the social planner's problem in \eqref{eq2:reg problem UB 3} and \eqref{eq2:social planner with DERs}.

It suffices to show now that $T^*_{\textsc{dec}}$ attains the upper bound in \eqref{eq2:reg problem UB 3}.
To that end, we use the independence sufficient condition $\omega \perp \lambda$.
We show that, under said condition, the expected total surplus $\overline{\sw}(T^*_{\textsc{dec}})$ matches the upper bound.
First note that the condition $\omega \perp \lambda$ allows to rewrite the upper bound in \eqref{eq2:reg problem UB 2} and \eqref{eq2:reg problem UB 3} as follows.
\begin{align*}
& \max_{\{q^i(\cdot),s^i(\cdot)\}_{i=1}^{M}} \ \Big\{ \ \overline{\sw} \ \big| \ s^i(\omega^i) \in \Umsc(\theta^i), \ i=1,\ldots,M \Big\} \nn\\
&= \sum_{i=1}^M \max_{q^i(\cdot)} \ \Embb_{\omega^i} \left[ S^i(q^i(\omega^i),\omega^i) - \Embb_{\lambda|\omega^i}\left[\lambda | \omega^i\right]^{\Top} q^i(\omega^i) \right] \nn\\
	&\qquad\qquad\qquad + \Embb_{\xi}\left[ \lambda^{\Top} r^i(\omega^i) \right] \nn\\
 	&\qquad\qquad\qquad +
\max_{s^i(\cdot)} \Big\{ \left. \Embb_{\xi}\left[ \lambda^{\Top} s^i(\omega^i) \right] \ \right| \  s^i(\omega^i) \in \Umsc(\theta^i) \Big\} \nn \displaybreak[0] \\
&= \sum_{i=1}^M \max_{q^i(\cdot)} \ \Embb_{\omega^i} \left[ S^i(q^i(\omega^i),\omega^i) - \overline{\lambda}^{\Top} q^i(\omega^i) \right] \quad + \quad \overline{\lambda}^{\Top} \overline{r}^i \\
 &\qquad\qquad \qquad +
\max_{s^i(\cdot)} \Big\{ \left. \overline{\lambda}^{\Top} \overline{s}^i \ \right| \  s^i(\omega^i) \in \Umsc(\theta^i) \Big\}  \displaybreak[0]  \\
&= \sum_{i=1}^M \Embb_{\omega^i} \left[ S^i(D^i(\overline{\lambda},\omega^i),\omega^i) - \overline{\lambda}^{\Top} D^i(\overline{\lambda},\omega^i) \right]  \\
	&\qquad\qquad\qquad + \overline{\lambda}^{\Top} \overline{r}^i  \quad + \quad \overline{\lambda}^{\Top} s^{*}(\overline{\lambda},\theta^i),
\end{align*}
where the last equality follows from the definition of the demand function $D^i(\pi,\omega^i)$ and the simplification of storage operation policy under deterministic prices.

The result follows since the tariff $T^*_{\textsc{dec}}$ induces the same expected total surplus if $\omega \perp \lambda$, \ie
\begin{align*} 
\overline{\sw}(T^*_{\textsc{dec}}) 
& = \Embb_{\xi} \bigg[ \sum_{i=1}^M S^i(D^{i}(\pi^*_{\textsc{dec}},\omega^i),\omega^i) - \lambda^{\Top} D^{i}(\pi^*_{\textsc{dec}},\omega^i)\bigg]  \\
& \qquad\qquad\qquad + \Embb \bigg[ \sum_{i=1}^M \lambda^{\Top} (r^i(\omega^i) + s^{*}(\pi^*_{\textsc{dec}},\theta^i))\bigg] \\
& = \sum_{i=1}^M \Embb_{\omega^i} \bigg[ S^i(D^{i}(\overline{\lambda},\omega^i),\omega^i) - \overline{\lambda}^{\Top} D^{i}(\overline{\lambda},\omega^i)\bigg] \\
& \qquad\qquad\qquad\qquad\qquad + \overline{\lambda}^{\Top} (\overline{r}^i + s^{*}(\overline{\lambda},\theta^i)).
\qedadhoc
\end{align*}

\subsubsection*{Proof of Theorem \ref{thm:centralized}} \label{proof:thm:centralized}

To solve problem \eqref{eq2:reg problem} over affine tariffs of the form $T(d)=A+\pi^{\Top}d$, we first need expressions for $\overline{\cs}(T)$ and $\overline{\rp}(T)$ in terms of $(\pi, A)$, considering the retailer-integrated DERs.

On the one hand, since this case does not consider customer-integrated DERs, the customers derive an expected surplus that remains unchanged, \ie $\overline{\cs}^i(T) = \overline{\cs}^i_0(T)$.
On the other hand, \eqref{eq2:retailer separation} characterizes the expected retailer surplus induced by $T$ considering the retailer-integrated DERs.

According to \eqref{eq2:retailer separation}, $\overline{\rp}(T)$ depends on the decision variables $(\pi,A)$ only through $\overline{\rp}_0(T)$.
Hence, the regulator's problem \eqref{eq2:reg problem} at hand, with parameter $F$, is equivalent to that of the regulator without any DERs and parameter $\tilde{F} := F-V(\overline{\lambda},\theta^o) - \Embb[ \lambda^{\Top} r^o(\xi)]$, \ie
\begin{align} \label{eq2:reg problem centralized}
\max_{T(\cdot)} \ \overline{\cs}_0(T) \quad \text{s.t.} \quad \overline{\rp}_0(T)=\tilde{F}.
\end{align}
Theorem 1 in \cite{MunozTong16partIarxiv} for the case without DERs characterizes the optimal two-part tariff for problem \eqref{eq2:reg problem centralized}.
Applying this result to problem \eqref{eq2:reg problem centralized} yields the desired result, \ie
$$
\pi^*_{\textsc{cen}} = \overline{\lambda} + \Embb[\nabla_{\pi}D(\pi^*_{\textsc{cen}},\omega)]^{-1} \Embb[\nabla_{\pi}D(\pi^*_{\textsc{cen}},\omega) (\lambda-\overline{\lambda})]
$$
and
\begin{align*}
A^*_{\textsc{cen}} &= \mbox{$\frac{1}{M}$} \left( \tilde{F} - \Embb \left[(\pi^*_{\textsc{cen}} - \lambda)^{\Top} D(\pi^*_{\textsc{cen}},\omega) \right] \right) \\
	&= \mbox{$\frac{1}{M}$} \left( F - \Embb \left[(\pi^*_{\textsc{cen}} - \lambda)^{\Top} D(\pi^*_{\textsc{cen}},\omega) \right] \right) \\
	& \qquad\qquad\qquad\qquad\qquad - \mbox{$\frac{1}{M}$} \left( V(\overline{\lambda},\theta^o) - \Embb[ \lambda^{\Top} r^o(\xi)] \right) \\
	&= A^* - \mbox{$\frac{1}{M}$} \left( V(\overline{\lambda},\theta^o) - \Embb[ \lambda^{\Top} r^o(\xi)] \right).
\qedadhoc
\end{align*}

\subsubsection*{Proof of Corollary \ref{cor:centralized}} \label{proof:cor:centralized}
Applying Cor. 2 in \cite{MunozTong16partIarxiv} to problem \eqref{eq2:reg problem centralized} implies that
\begin{align*}
\overline{\cs}(T^*_{\textsc{cen}}) &= \overline{\sw}^*_0 - \tilde{F} \\
	&= \overline{\sw}^*_0 - F + V(\overline{\lambda},\theta^o) + \Embb[ \lambda^{\Top} r^o(\xi)]
\end{align*}
where, according to \eqref{eq2:sw0 star}, $\overline{\sw}^*_0 \equiv \overline{\sw}_0(T^*)$ does not depend on $F$.
The result follows since $\overline{\rp}(T^*_{\textsc{cen}})=F$ further implies that
\begin{align*}
\overline{\sw}(T^*_{\textsc{cen}}) &= \overline{\cs}(T^*_{\textsc{cen}}) + \overline{\rp}(T^*_{\textsc{cen}}) \\
		&= \overline{\sw}^*_0 + V(\overline{\lambda},\theta^o) + \Embb[ \lambda^{\Top} r^o(\xi)].
\qedadhoc
\end{align*}

\end{document}